\documentclass[11pt]{article}
\usepackage{graphicx}
\usepackage{hyperref}
\usepackage{amsmath,amsthm,amssymb,enumerate}
\usepackage[normalem]{ulem} 
\usepackage{euscript,mathrsfs}
\usepackage[left=3cm,right=3cm,top=3.5cm,bottom=3.5cm]{geometry}
\usepackage{color}
\usepackage{bm}
\usepackage{bbm}
\catcode`\@=11 \@addtoreset{equation}{section}

\catcode`\@=12
\usepackage{color}
\allowdisplaybreaks

\newtheorem{Theorem}{Theorem}[section]
\newtheorem{Proposition}[Theorem]{Proposition}
\newtheorem{Lemma}[Theorem]{Lemma}
\newtheorem{Corollary}[Theorem]{Corollary}

\theoremstyle{definition}
\newtheorem{Definition}[Theorem]{Definition}

\newtheorem{Remark}[Theorem]{Remark}

\newcommand{\bTheorem}[1]{
\begin{Theorem} \label{T#1} }
\newcommand{\eT}{\end{Theorem}}

\newcommand{\bProposition}[1]{
\begin{Proposition} \label{P#1}}
\newcommand{\eP}{\end{Proposition}}

\newcommand{\bLemma}[1]{
\begin{Lemma} \label{L#1} }
\newcommand{\eL}{\end{Lemma}}

\newcommand{\bCorollary}[1]{
\begin{Corollary} \label{C#1} }
\newcommand{\eC}{\end{Corollary}}

\newcommand{\bRemark}[1]{
\begin{Remark} \label{R#1} }
\newcommand{\eR}{\end{Remark}}

\newcommand{\bDefinition}[1]{
\begin{Definition} \label{D#1} }
\newcommand{\eD}{\end{Definition}}

\newcommand{\dif}{\mathrm{d}}

\newcommand{\mr}{\mathbb{R}}
\newcommand{\prst}{\mathbb{P}}
\newcommand{\p}{\mathbb{P}}

\newcommand{\mn}{\mathbb{N}}

\newcommand{\mt}{\mathbb{R}^n}

\newcommand{\bq}{\mathbf q}

\newcommand{\tor}{\mathbb{R}^n}

\newcommand{\StoB}{\left(\Omega, \mathfrak{F},(\mathfrak{F}_t )_{t \geq 0},  \mathbb{P}\right)}

\newcommand{\bfu}{\mathbf{u}}
\newcommand{\bfv}{\mathbf{v}}
\newcommand{\bfq}{\mathbf{q}}

\newcommand{\bfa}{\mathbf{a}}

\newcommand{\ds}{\,\mathrm{d}\sigma}

\newcommand{\bFormula}[1]{
\begin{equation} \label{#1}}
\newcommand{\eF}{\end{equation}}

\newcommand{\vr}{\varrho}

\newcommand{\vu}{\vc{u}}
\newcommand{\vc}[1]{{\bf #1}}

\newcommand{\Div}{{\rm div}_x}
\newcommand{\Grad}{\nabla_x}

\newcommand{\tn}[1]{\mathbb{#1}}
\newcommand{\dx}{\,{\rm d} {x}}

\newcommand{\dt}{\,{\rm d} t }

\newcommand{\D}{{\rm d}}
\newcommand{\ep}{\varepsilon}

\newcommand{\R}{\mathbb{R}}

\newcommand{\E}{\mathbb{E}}

\definecolor{Cgrey}{rgb}{0.85,0.85,0.85}
\definecolor{Cblue}{rgb}{0.50,0.85,0.85}
\definecolor{Cred}{rgb}{1,0,0}
\definecolor{fancy}{rgb}{0.10,0.85,0.10}

\newcommand\Cbox[2]{%
    \newbox\contentbox%
    \newbox\bkgdbox%
    \setbox\contentbox\hbox to \hsize{%
        \vtop{
            \kern\columnsep
            \hbox to \hsize{%
                \kern\columnsep%
                \advance\hsize by -2\columnsep%
                \setlength{\textwidth}{\hsize}%
                \vbox{
                    \parskip=\baselineskip
                    \parindent=0bp
                    #2
                }%
                \kern\columnsep%
            }%
            \kern\columnsep%
        }%
    }%
    \setbox\bkgdbox\vbox{
        \color{#1}
        \hrule width  \wd\contentbox %
               height \ht\contentbox %
               depth  \dp\contentbox
        \color{black}
    }%
    \wd\bkgdbox=0bp%
    \vbox{\hbox to \hsize{\box\bkgdbox\box\contentbox}}%
    \vskip\baselineskip%
}

\date{}

\begin{document}


\title{Stochastic compressible Euler equations and inviscid limits}

\author{Dominic Breit 
\and Prince Romeo Mensah}

\date{\today}

\maketitle

\centerline{Department of Mathematics, Heriot-Watt University}

\centerline{Riccarton Edinburgh EH14 4AS, UK}

%
%
%
%
%

\bigskip


\begin{abstract}

We prove the existence of a unique local strong solution to the stochastic compressible Euler system with nonlinear multiplicative noise. This solution
exists up to a positive stopping time and is strong in both the PDE and probabilistic sense. Based on this existence result, we study the inviscid limit of the stochastic compressible Navier--Stokes system. As the viscosity tends to zero, any sequence of finite energy weak martingale solutions converges to the compressible Euler system.

%

\end{abstract}



\section{Introduction}
\label{P}

We consider a stochastic variant of the \emph{compressible barotropic
Euler system} describing the time evolution of the mass density $\vr$ and the bulk velocity $\vu$ of a fluid driven by
a nonlinear multiplicative noise. The system of equations reads
\begin{align} \label{P1}
\D \vr + \Div (\vr \vu) \dt &= 0,\\ \label{P2}
\D (\vr \vu) + \left[ \Div (\vr \vu \otimes \vu) + a \Grad \vr^\gamma \right]  \dt  &=  \mathbb{G} (\vr, \vr \vu) \D W.
\end{align}
Here $\gamma>1$ denotes the adiabatic exponent, $a>0$ is the squared reciprocal of the Mach number (the ratio between average velocity and speed of sound).
The driving process $W$ is a cylindrical Wiener process defined on some probability space $(\Omega,\mathfrak{F},\p)$ and the coefficient $\mathbb{G}$ is generally nonlinear and satisfies suitable growth assumptions, see Section \ref{E} for the precise set-up. In order to eliminate the well-known difficulties related to the behaviour of fluid flows near the boundary of the underlying domain but still consider a physically meaningful situation, we study \eqref{P1}--\eqref{P2} on the whole space $\R^n$. We complement \eqref{P1}--\eqref{P2} with the far field condition
\begin{align}\label{farfield}
\varrho(x)\rightarrow \overline\varrho,\quad \bfu(x)\rightarrow 0,\quad |x|\rightarrow\infty,
\end{align}
for some $\overline\varrho>0$. 
The initial conditions are random variables
\begin{equation} \label{P4}
\vr(0,\cdot) = \vr_0, \ \vu(0, \cdot) = \vu_0,
\end{equation}
with sufficient spatial regularity specified later. The main interest is the three-dimensional case, but $n=1,2$ are also included in our theory (and obviously higher but non-physical dimensions). We remark that in contrast to the incompressible system,  the one-dimensional situation makes sense for \eqref{P1}--\eqref{P2}. 
In fact, martingale solutions to \eqref{P1}--\eqref{P2} for $n=1$ are studied in \cite{BV13}.\\\

Our main result concerning the system  \eqref{P1}--\eqref{P4} is the existence
of a unique \emph{maximal strong pathwise solution}. This solution is strong in the analytical sense (i.e., equations \eqref{P1}--\eqref{P2} are satisfied pointwise) and strong in the probabilistic sense (i.e., it is defined on a given probability space). It exists up to the hypothetical blow-up of the $W^{1,\infty}$-norm of the velocity $\vu$. The precise formulation is given in Definition \ref{def:maxsol}. The existence of a \emph{maximal strong pathwise solution}, being defined on a maximal (random) time interval, follows from an extension of a \emph{local strong pathwise solution}, 
see Definition \ref{def:strsol}, which lives up to a suitable stopping time. The main statement can be found in Theorem \ref{thm:main}. Corresponding results
in the deterministic case are classical and we refer to \cite{Ag} and \cite{hugo}.\\
 As in the incompressible case, global existence and uniqueness is a famous open problem. The presence of noise does not seem to change the situation. As solutions to nonlinear hyperbolic systems are known to develop singularities in finite time, the question about  global well-posedness in the class of weak solutions has been analysed extensively. This is based on the method of convex integration which has been developed in the context of fluid mechanics by De Lellis and Sz\'ekelyhidi \cite{DelSze}. The non-uniqueness of global-in-time weak solutions to \eqref{P1}--\eqref{P2} has recently  been shown in \cite{BrFeHo2017}
proceeding similar result in the deterministic case, cf. \cite{Fei2016}.\\
In contrast to the compressible system \eqref{P1}--\eqref{P2}, its incompressible counterpart has been studied extensively. There are numerous results about the two-dimensional situation, see \cite{Bes,BeFl,BrPe,CaCu,Kim2}. First results in three dimensions (treated on the whole space) can be found in \cite{MiVa,Kim}. Similar to
our main theorem, the existence of a unique local strong solutions is shown, however, only additive noise is allowed.
The general three-dimensional case (with slip boundary conditions in a bounded domain and with nonlinear multiplicative noise) has finally been studied recently in \cite{GHVic}.\\
A main idea in our existence proof is to rewrite \eqref{P1}--\eqref{P2} as a symmetric hyperbolic system by formally dividing \eqref{P2} by $\varrho$ similarly to \cite{BrFeHo2016} . In order to make the general framework from \cite{JUKim} for these systems available, we cut the noise in the critical range (that is, if $\varrho$ is large or close to zero). The main tool in the limit procedure is an abstract Cauchy lemma from \cite{GHZi}, see Lemma \ref{lem:GH}. We remark that the method from \cite{BrFeHo2016}, used in the analysis of the Navier--Stokes,
can to a certain extend be applied to \eqref{P1}--\eqref{P2} at least if periodic boundary conditions are considered. It does not, however, yield continuity of density and velocity in time -- an advantage of the approach in the present paper.
\\\

In our second main result we are concerned with the relationship between the Navier--Stokes and Euler equations. Viscous compressible fluids subject to stochastic forcing
can be described by the Navier--Stokes system
\begin{align} \label{P1NS}
\D \vr + \Div (\vr \vu) \dt &= 0,\\
 \label{P2NS}
\D (\vr \vu) + \left[ \Div (\vr \vu \otimes \vu) + a \Grad \vr^\gamma \right]  \dt  &= \Div \mathbb{S} (\Grad \vu) \ \dt + \mathbb{G} (\vr, \vr \vu) \D W.
\end{align}
Here
$\mathbb{S}(\Grad \vu)$ is the viscous stress tensor for which we assume Newton's rheological law
\begin{equation} \label{P3NS}
\mathbb{S}(\Grad \vu) = \nu \left( \Grad \vu + \Grad^t \vu - \frac{2}{n} \Div \vu \mathbb{I} \right) + \lambda \Div \vu \mathbb{I}, \qquad \nu > 0, \ \lambda \geq 0.
\end{equation}
The study of the  system  \eqref{P1NS}--\eqref{P3NS} was first initiated in \cite{BrHo} where the global-in-time existence of \emph{finite energy weak martingale solutions} is shown. These solutions are weak in the analytical sense (derivatives only exists in the sense of distributions) and weak in the probabilistic sense (the probability space is an integral part of the solution) as well. Moreover, the time-evolution of the energy can be controlled in terms of its initial state. The results from \cite{BrHo} -- limited to periodic boundary conditions -- have been extended to the whole space in \cite{Romeo1}. So, a comparison between \eqref{P1}--\eqref{P2} and 
\eqref{P1NS}--\eqref{P3NS} is possible. Our main result, stated in Theorem \ref{thm:invisicdLimit}, shows that any sequence of \emph{finite energy weak martingale solutions} to \eqref{P1NS}--\eqref{P3NS} converges locally in time to the unique strong solution of \eqref{P1}--\eqref{P2} as $\lambda,\nu\rightarrow0$. 
A similar strategy has been employed in \cite{BrFeHo2015A} in order to study the inviscid-incompressible limit (where in addition, $a=\frac{1}{\varepsilon^2}$ with $\varepsilon\rightarrow0$ is considered, where the limit system is the incompressible Euler system).
A major difference to \cite{BrFeHo2015A} is the generality of the noise coefficients we can consider now.
Due to the incompressibility constraint on the target system only linear noise can be considered  in \cite{BrFeHo2015A}.
In contrast to this, in the compressible case we can allow the full generality for the noise for which the existence theory applies.\\ The main tool in our proof is the relative energy inequality from \cite{BrFeHo2015A}. It allows to compare a \emph{finite energy weak martingale solutions} to \eqref{P1NS}--\eqref{P3NS}
with a set of smooth comparison functions -- 
in this case, the local solution to \eqref{P1}--\eqref{P2}. The concept of relative energy inequality has a long history starting with the pioneering work of Dafermos \cite{Daf4}.
In the context of compressible Navier--Stokes equation, it has been introduced in \cite{FeJaNo}
and has also been  used to study the inviscid limit of the compresible Navier--Stokes system in the deterministic case, cf. \cite{Sueur}.

\section{Preliminaries and main result} \label{E}

We start by introducing some notations and  basic facts used in the text. To begin, we fix
an arbitrary large time horizon $T>0$.

\subsection{Analytic framework}
We will define the Sobolev space $W^{s,2}(\mathbb{R}^n)$ for $s\in\mathbb{R}$ as the set of tempered distributions for which the norm
\begin{align}
\Vert  v\Vert_{W^{s,2}(\mathbb{R}^n)} = \bigg(\int_{\mathbb{R}^n} \big(1+\vert \xi\vert^2  \big)^s\vert  \hat{v}(\xi)\vert^2\, \mathrm{d}\xi
  \bigg)^\frac{1}{2}
\end{align}
defined in frequency space is finite. Here $\hat{v}$ denotes the Fourier transform of $v$. To shorten notation, we will write $\Vert  \cdot\Vert_{s,2}$ for $\Vert \cdot\Vert_{W^{s,2}(\mathbb{R}^n)}$.
The following estimates are standard in the Moser-type calculus and can be found e.g. in Majda \cite[Proposition 2.1]{Majd}.

\begin{enumerate}

\item For $u,v \in W^{s,2} \cap L^\infty(\mt)$ and $|\alpha| \leq s$ we have
\begin{equation}\label{E5}
\left\| \partial^\alpha_x (u v) \right\|_{2} \leq c_s \left( \| u \|_{\infty} \| \nabla^s_x v \|_{2} +
\| v \|_{\infty} \| \nabla^s_x u \|_{2} \right).
\end{equation}

\item For $u \in W^{s,2}(\mt)$, $\nabla_x u \in L^\infty(\mt)$, $v \in W^{s-1,2} \cap L^\infty (\mt)$ and $|\alpha| \leq s$ we have
\begin{equation}\label{E6}
\left\| \partial^\alpha_x (uv) - u \partial^\alpha_x v \right\|_{2} \leq c_s
\left( \| \nabla_x u \|_{\infty} \| \nabla^{s-1}_x v \|_{2} +
\| v \|_{\infty} \| \nabla^s_x u \|_{2} \right).
\end{equation}

\item Let $u \in W^{s,2} \cap C(\mt)$ and let $F$ be an $s$-times continuously differentiable function on an open neighborhood of the compact set $G =
{\rm range}[u]$. Then we have for all 
$1 \leq |\alpha| \leq s$,
\begin{equation} \label{E7}
\left\| \partial^\alpha_x F(u) \right\|_{L^2(\mt)} \leq c_s \| \partial_u F \|_{C^{s-1}(G)} \| u \|^{|\alpha| - 1}_{L^\infty(\mt)} \|
\partial^\alpha_x u \|_{L^2(\mt)}.
\end{equation}
\color{black}
\end{enumerate}

\subsection{Stochastic framework}

The driving process $W$ is a cylindrical Wiener process on a separable Hilbert space $\mathfrak{U}$ defined on some stochastic basis $\StoB$ with a complete, right-continuous filtration. More specifically, $W$ is given by a formal expansion
\[
W(t)=\sum_{k\geq 1} e_k \beta_k(t),
\]
where $\{ \beta_k \}_{k \geq 1}$ is a family of mutually independent real-valued Brownian motions
with respect to $\StoB$ and $\{e_k\}_{k\geq 1}$ is an orthonormal basis of  $\mathfrak{U}$.
To give the precise definition of the diffusion coefficient $\mathbb{G}$, consider $\rho\in L^2(\mt)$, $\rho\geq0$, $\bfq\in L^2(\mt)$ and define it as follows
$$\mathbb{G}(\rho,\bq)e_k=\mathbf{G}_k(\cdot,\rho(\cdot),\bq(\cdot)).$$
We suppose that
the coefficients $\mathbf{G}_{k}:\mt\times [0,\infty) \times\mr^n\rightarrow\mr^n$ are $C^s$-functions that satisfy uniformly in $x\in\mt$ 
\begin{equation}\label{FG1}
\vc{G}_k (\cdot, 0 , 0) = 0,
\end{equation}
\begin{equation}
|\nabla^l \vc{G}_k (\cdot, \cdot, \cdot) | \leq \alpha_k, \quad \sum_{k \geq 1} \alpha_k  < \infty \quad \mbox{for all}\ l\in\{1,...,s\},
\label{FG2}
\end{equation}
with $s\in\mn$ specified below. Finally, we assume that the $\vc{G}_k$s are compactly supported, i.e. there is $\mathcal K\Subset\R^n$ such that
\begin{align}
\label{FG3}\mathrm{spt}(\vc{G}_k)\Subset\mathcal K\quad\text{for all}\quad k\in \mathbb N.
\end{align}
This is also assumed in the Navier--Stokes case in view of the far field condition, cf. \cite{Romeo1}.
A typical example we have in mind is
\begin{align}\label{eq:model}
\mathbf{G}_k(x,\rho,\bq)=\bfa_k(x)\rho+\mathbb A_k(x) \bfq,
\end{align}
where $\bfa_k:\mt\rightarrow\R^n$ and $\mathbb A_k:\mt\rightarrow\R^{n\times n}$ are smooth functions which are compactly supported. However, our analysis applies to  general nonlinear  coefficients $\mathbf{G}_k$.
%
%
%
%
%
%

\medskip

Observe that if $\vr$, $\vc{q}$ are $(\mathfrak{F}_t)$-progressively measurable $L^2(\mt)$-valued processes such that
\[
\vr \in L^2 \Big( \Omega\times[0,T]; L^2(\mt) \Big), \ \vc{q} \in L^2 \Big(\Omega\times[0,T]; L^2(\mt; \mr^n) \Big),
\]
and $\mathbb{G}$ satisfies (\ref{FG1}) and (\ref{FG2}), then
the stochastic integral
\[
\int_0^t \mathbb{G}(\vr, \vr \vu) \ {\rm d} W = \sum_{k \geq 1}\int_0^t \vc{G}_k (\cdot, \vr, \vr \vu) \ {\rm d} \beta_k
\]
is a well-defined  $(\mathfrak{F}_t)$-martingale ranging in $L^2(\mt; \mr^n)$.
Finally, we define an auxiliary space $\mathfrak{U}_0\supset\mathfrak{U}$ via
$$\mathfrak{U}_0=\bigg\{v=\sum_{k\geq1}\alpha_k e_k;\;\sum_{k\geq1}\frac{\alpha_k^2}{k^2}<\infty\bigg\},$$
endowed with the norm
$$\|v\|^2_{\mathfrak{U}_0}=\sum_{k\geq1}\frac{\alpha_k^2}{k^2},\quad v=\sum_{k\geq1}\alpha_k e_k.$$
Note that the embedding $\mathfrak{U}\hookrightarrow\mathfrak{U}_0$ is Hilbert-Schmidt. Moreover, trajectories of $W$ are $\prst$-a.s. in $C([0,T];\mathfrak{U}_0)$, cf. \cite{daprato}.\\
Let us complete this section by presenting a technical tool to pass to the limit in sequences of local strong solutions (which typically exists only up to a stopping time). It originates from \cite[Lemma 5.1]{GHZi} in an abstract setting. For the present form, we refer to \cite[Lemma 7.1]{GHVic}.
\begin{Lemma}
\label{lem:GH}
Let $\StoB$ be a stochastic basis with a complete right-continuous filtration.
Let $(u_R)$ be a sequence of $W^{s,2}(\mt)$-valued continuous stochastic processes adapted to $(\mathfrak F_t)$. For $N>1$ and $T>0$, we define the sequence of stopping times
\begin{align*}
\mathfrak t_{R,N}=\inf\big\{t\geq0:\,\,\|u_R(t)\|_{s,2}\geq N+\|u_R(0)\|_{s,2}\big\}\wedge T.
\end{align*}
 Assume that we have
\begin{align}\label{GH1}
&\lim_{R\rightarrow\infty}\sup_{L\geq R}\E\sup_{0\leq t\leq \mathfrak t_{R,N}\wedge \mathfrak t_{L,N}}\|u_R(t)-u_L(t)\|_{s,2}=0,\\
\label{GH2}&\lim_{\delta\rightarrow0}\sup_{ R}\p\bigg\{\sup_{0\leq t\leq \mathfrak t_{R,N}\wedge \delta}\|u_R(t)\|_{s,2}>\|u_R(0)\|_{s,2}+N-1\bigg\}=0.
\end{align}
Then there is a stopping time $\mathfrak t$ such that $\mathfrak t>0$ $\p$-a.s. and
a $W^{s,2}(\mt)$-valued $(\mathfrak{F}_t)$-progressively measurable process $u$ satisfying
$$ u(\cdot\wedge \mathfrak{t}) \in   C([0,T]; W^{s,2}(\tor))\quad \mathbb{P}\text{-a.s.}$$
such that
\begin{align}\label{GH3}
\sup_{0\leq t\leq \mathfrak t}\|u_R(t)-u(t)\|_{s,2}=0
\end{align}
 $\p$-a.s. as $R\rightarrow\infty$ (at least for a subsequence).
\end{Lemma}

\subsection{Compressible Euler equations}

Let us first introduce the notion of local strong pathwise solution. Such a solution is strong in both the PDE and probabilistic sense but possibly exists only locally in time. To be more precise, system \eqref{P1}--\eqref{P2} will be satisfied pointwise (not only in the sense of distributions) on the given stochastic basis associated to the cylindrical Wiener process $W$.

\begin{Definition}[Local strong pathwise solution] \label{def:strsol}

Let $\StoB$ be a stochastic basis with a complete right-continuous filtration. Let ${W}$ be an $(\mathfrak{F}_t) $-cylindrical Wiener process and let $\overline\varrho>0$.
Let
$(\varrho_0,\bfu_0)$ be a $\overline\vr+W^{s,2}(\mt)\times W^{s,2}(\mt)$-valued $\mathfrak{F}_0$-measurable random variable, and let $\mathbb{G}$ satisfy \eqref{FG1} and
\eqref{FG2} for some $s\in\mathbb N$.
A triplet
$(\varrho,\vu,\mathfrak{t})$ is called a local strong pathwise solution to the system \eqref{P1}--\eqref{P4} provided
\begin{enumerate}[(a)]
\item $\mathfrak{t}$ is an a.s. strictly positive  $(\mathfrak{F}_t)$-stopping time;
\item the density $\varrho$ is a $\overline\varrho+W^{s,2}(\mt)$-valued $(\mathfrak{F}_t)$-progressively measurable process satisfying
$$\varrho(\cdot\wedge \mathfrak{t})  > 0,\ \varrho(\cdot\wedge \mathfrak{t}) \in C([0,T]; \overline\varrho+W^{s,2}(\tor)) \quad \mathbb{P}\text{-a.s.};$$
\item the velocity $\vu$ is a $W^{s,2}(\mt)$-valued $(\mathfrak{F}_t)$-progressively measurable process satisfying
$$ \vu(\cdot\wedge \mathfrak{t}) \in   C([0,T]; W^{s,2}(\tor))\quad \mathbb{P}\text{-a.s.};$$
\item  there holds $\prst$-a.s.
\[
\begin{split}
\varrho (t\wedge \mathfrak{t}) &= \varrho_0 -  \int_0^{t \wedge \mathfrak{t}} \Div(\varrho\vu ) \ \dif s, \\
(\varrho \vu) (t \wedge \mathfrak{t})  &= \varrho_0 \vu_0 - \int_0^{t \wedge \mathfrak{t}} \Div (\varrho\vu \otimes\vu ) \ \dif s 
- \int_0^{t \wedge \mathfrak{t}}a\Grad \varrho^\gamma \ \dif s + \int_0^{t \wedge \mathfrak{t}} {\tn{G}}(\varrho,\varrho\vu ) \ \D W,
\end{split}
\]
for all $t\in[0,T]$.
\end{enumerate}
\end{Definition}

In the above definition, we have tacitly assumed that $s$ is large enough in order to provide sufficient regularity for the strong solutions.
Classical solutions require spatial derivatives of $\vu$ and $\vr$ to be continuous $\prst$-a.s. This motivates the following definition.

\begin{Definition}[Maximal strong pathwise solution]\label{def:maxsol}
Fix a stochastic basis with a cylindrical Wiener process and an initial condition as in Definition \ref{def:strsol}. A quadruplet $$(\varrho,\vu,(\mathfrak{t}_R)_{R\in\mn},\mathfrak{t})$$ is a maximal strong pathwise solution to system \eqref{P1}--\eqref{P4} provided

\begin{enumerate}[(a)]
\item $\mathfrak{t}$ is an a.s. strictly positive $(\mathfrak{F}_t)$-stopping time;
\item $(\mathfrak{t}_R)_{R\in\mn}$ is an increasing sequence of $(\mathfrak{F}_t)$-stopping times such that
$\mathfrak{t}_R<\mathfrak{t}$ on the set $[\mathfrak{t}<T]$,
$\lim_{R\to\infty}\mathfrak{t}_R=\mathfrak t$ a.s. and
\begin{equation}\label{eq:blowup}
\sup_{t\in[0,\mathfrak{t}_R]}\|\vu(t)\|_{1,\infty}\geq R\quad \text{on}\quad [\mathfrak{t}<T] ;
\end{equation}
\item each triplet $(\varrho,\vu,\mathfrak{t}_R)$, $R\in\mn$,  is a local strong pathwise solution in the sense  of Definition \ref{def:strsol}.
\end{enumerate}
\end{Definition}

The notion of a maximal pathwise solution has already appeared in the literature in the context of various SPDE or SDE models, see for instance \cite{BMS,Elw, jac,MikRoz}.
%
%
Finally, we have all in hand to formulate our main result.

\begin{Theorem}\label{thm:main}
Let $s\in\mn$ satisfy $s>\frac{n}{2} + 2$ and let $\overline\varrho>0$.
Let the coefficients
$\mathbf{G}_k$ satisfy hypotheses \eqref{FG1}, \eqref{FG2} and let $(\varrho_0,\bfu_0)$ be an $\mathfrak{F}_0$-measurable, $\overline\vr+W^{s,2}(\mt)\times W^{s,2}(\mt)$-valued random variable such that $\varrho_0>0$ $\p$-a.s.
Then
there exists a unique maximal strong pathwise solution $(\varrho,\vu,(\mathfrak{t}_R)_{R\in\mn},\mathfrak{t})$ to problem \eqref{P1}--\eqref{P4}
in the sense of Definition \ref{def:maxsol} with the initial condition $(\varrho_0,\vu_0)$.
\end{Theorem}

\begin{Remark}
Starting with the pioneering work in \cite{DelSze3}, several counterexamples have been developed showing that the (deterministic) compressible Euler system is desperately ill-posed. Even if the initial data is smooth, the global existence and uniqueness of solutions can fail.
Similar result for the stochastic compressible Euler system have been achieved recently in \cite{BrFeHo2017}. The existence of global strong solutions (i.e. the stopping time $\mathfrak t$ in Definition \ref{def:strsol} reaches $T$) is not expected.
\end{Remark}

\subsection{Compressible Navier--Stokes equations}
In this section, we present the concept of finite energy weak martingale solutions
to \eqref{P1NS}--\eqref{P3NS}. It has been introduced in \cite{BrHo}
and improved in \cite{BrFeHo2015A} and \cite{BrFeHobook}. Both papers complement \eqref{P1NS}--\eqref{P3NS} with periodic boundary conditions. A corresponding version on the whole space can be found in \cite{Romeo1}. These solutions are weak in the analytical sense (derivatives only exists in the sense of distributions) and weak in the probabilistic sense (the probability space is an integral part of the solution) as well. Moreover, the time-evolution of the energy can be controlled in terms of its initial state. They exists globally in time.

\begin{Definition}[Finite energy weak martingale solution]
\label{def:weakMartin}
Let $\overline\rho>0$.
Let $\Lambda$ be a Borel probability measure on $L^\gamma_{\mathrm loc}(\mathbb{R}^n)\times L_{\mathrm loc}^{2\gamma/\gamma+1}(\mathbb{R}^n)$. Then $\big[ \big( \Omega,\mathfrak{F}, (\mathfrak{F}_t)_{t\geq0},\mathbb{P} \big);\rho,\mathbf{v}, W \big]$ is a finite energy weak martingale solution of \eqref{P1NS}--\eqref{P3NS} if
\begin{enumerate}[(a)]
\item $ \big( \Omega,\mathfrak{F}, (\mathfrak{F}_t)_{t\geq0},\mathbb{P} \big)$ is a stochastic basis with a complete right-continuous filtration,
\item $W$ is a $(\mathfrak{F}_t)$-cylindrical Wiener process,
\item the density $\rho$ satisfies $\rho\geq 0$, $t\mapsto \langle \rho(t, \cdot),\phi\rangle\in C[0,T]$ for any $\phi\in C_c^\infty(\mathbb{R}^n)$ $\mathbb{P}$-a.s., the function $t\mapsto \langle \rho(t, \cdot),\phi\rangle$ is progressively measurable and 
\begin{align*}
\mathbb{E}\, \bigg[ \sup_{t\in(0,T)}\Vert  \rho(t,\cdot)\Vert_{L^\gamma(K)}^p\bigg]<\infty 
\end{align*}
for all $1\leq p<\infty$ and all $K\Subset\R^n$,
\item the velocity field $\mathbf{v}$ is $(\mathfrak{F}_t)$-adapted and 
\begin{align*}
\mathbb{E}\, \bigg[ \int_0^T\Vert  \mathbf{v}\Vert^2_{W^{1,2}(K)}\,\mathrm{d}t\bigg]^p<\infty 
\end{align*}
for all $1\leq p<\infty$ and all $K\Subset \R^n$,
\item the momentum $\rho\mathbf{v}$ satisfies  $t\mapsto \langle \rho\mathbf{v},\bm\varphi\rangle\in C[0,T]$ for any $\bm{\varphi}\in C_c^\infty(\mathbb{R}^n)$ $\mathbb{P}$-a.s., the function $t\mapsto \langle \rho\mathbf{v},\bm{\phi}\rangle$ is progressively measurable and 
\begin{align*}
\mathbb{E}\, \bigg[ \sup_{t\in(0,T)}\Vert  \rho\mathbf{v}\Vert_{L^\frac{2\gamma}{\gamma+1}(K)}^p\bigg]<\infty 
\end{align*}
for all $1\leq p<\infty$ and all $K\Subset \R^n$,
\item $\Lambda=\mathbb{P}\circ \big[ \rho(0),\rho\mathbf{v}(0) \big]^{-1}$,
\item for all $\phi\in C^\infty_c(\mathbb{R}^n)$ and $\bm{\varphi}\in C_c^\infty(\mathbb{R}^n)$ we have
\begin{equation}\label{eq:energy}
\begin{aligned}
\langle \rho(t), \phi\rangle &= \langle\rho(0) , \phi\rangle - \int_0^{t}\langle\rho\mathbf{v}, \nabla_x \phi\rangle\mathrm{d}s 
\\
\langle \rho\mathbf{v}(t), \bm\varphi\rangle &= \langle(\rho\mathbf{v})(0), \bm{\varphi}\rangle - \int_0^{t}\langle \rho\mathbf{v}\otimes\mathbf{v}, \nabla_x  \bm{\varphi}\rangle\,\mathrm{d}s
+ \int_0^{t} \langle\mathbb{S}(\nabla_x \mathbf{v})\,, \nabla_x  \bm{\varphi}\rangle\mathrm{d}s 
\\&
- \int_0^{t}\langle a \rho^\gamma , \mathrm{div}_x\,\bm{\varphi}\rangle\,\mathrm{d}s
+\int_0^{t}\langle\mathbb{G}(\rho,\rho\mathbf{v}), \bm{\varphi}\rangle\,\mathrm{d}W
\end{aligned}
\end{equation}
$\mathbb{P}$-a.s. for all $t\in[0,T]$,
\item the energy inequality\index{energy inequality}
\begin{align}\label{EI2''}
\begin{aligned}
&
\int_{\R^n} \Big[ \frac{1}{2} \varrho | {\bf v} |^2  + H(\rho,\overline\rho) \Big] \dx 
+ \int_0^t \int_{\R^n}  \mathbb{S} (\nabla_x  {\bf v}): \nabla_x  {\bf v} \dx\,  {\rm d}s \\
&\leq \int_{\R^n} \Big[ \frac{1}{2}\frac{|\rho\bfv(0)|^2}{\rho(0)}  + H(\rho(0),\overline\rho) \Big] \dx
+\sum_{k=1}^\infty\int_0^t\bigg(\int_{\R^n}\mathbf{G}_k (\rho, \rho {\bf v})\cdot{\bf v}\dx\bigg){\rm d} W_k\\
&+ \frac{1}{2}\sum_{k = 1}^{\infty}  \int_0^t
\int_{\R^n} \rho^{-1}| \mathbf{G}_k (\rho, \rho {\bf v}) |^2 \, {\rm d}s
\end{aligned}
\end{align}
holds for a.e. $t\in[0,T]$ $\mathbb P$-a.s. Here $M_E$ is a real-valued square integrable martingale 
and $H(\rho,\overline\rho)$ is the pressure potential given by
\begin{align}
H(\rho,\overline \rho)=\frac{a}{\gamma-1}\Big(\rho^\gamma-\gamma\overline\rho^{\gamma-1}(\rho-\overline\rho)-\overline\rho^\gamma\Big).
\end{align}
\end{enumerate}
\end{Definition}
The following existence theorem is shown in \cite[Theorem 1]{Romeo1} (see also \cite[Ch. 3]{Romeo2} for a more detailed proof).
\begin{Theorem}
\label{thm:Romeo1}
Let $\overline\rho>0$, $\gamma>\frac{n}{2}$ and assume that $\Lambda$ is a Borel probability measure on $L^\gamma_{\mathrm loc}(\mathbb{R}^n)\times L_{\mathrm loc}^{2\gamma/\gamma+1}(\mathbb{R}^n)$ such that
\begin{align*}
\Lambda\Bigg\{(\rho,\mathbf{q})\in L^\gamma_{\mathrm loc}(\mathbb{R}^n)\times L_{\mathrm loc}^{2\gamma/\gamma+1}(\mathbb{R}^n) \, :\, 0<M_1\leq \rho\leq M_2,\, \mathbf{q}\vert_{\{\rho=0\}}=0,   \Bigg\}=1  
\end{align*}
with constants $0<M_1<M_2$. Furthermore, assume that the following moment estimate
\begin{align}\label{initial}
\int_{L^\gamma_x\times  L^{2\gamma/\gamma+1}_x} \bigg\Vert\frac{1}{2}\frac{\vert  \mathbf{q} \vert^2}{\rho}+ H(\rho,\overline\rho) \bigg\Vert^p_{L^1_x}\, \mathrm{d}\Lambda(\rho,\mathbf{q})<\infty
\end{align}
holds for all $1\leq p<\infty$. Finally, assume that \eqref{FG1}, \eqref{FG2} with $s=1$ and \eqref{FG3}  hold. Then there exists a finite energy weak martingale solution of \eqref{P1NS}--\eqref{P3NS} in the sense of Definition \ref{def:weakMartin} with initial law $\Lambda$.
\end{Theorem}
\begin{Remark}\label{rem:new}
Due to the assumptions on the initial law as well as \eqref{FG1}, \eqref{FG2} with $s=1$ and \eqref{FG3} the energy inequality \eqref{EI2''} implies the following moment estimates
\begin{align*}
\E\bigg[\int_0^T\int_{\R^n}|\nabla {\bf v}|^2\dx \dt\bigg]^p<\infty,\quad
\E\bigg[\int_0^T\int_{\R^n}\rho|{\bf v}|^2\dx \dt\bigg]^p<\infty,\\
\E\bigg[\int_0^T\int_{\R^n}\Big(|\rho-\overline\rho|^2\mathbb I_{|\rho-\overline\rho|\leq  1}+|\rho-\overline\rho|^\gamma\mathbb I_{|\rho-\overline\rho|> 1}\Big)\dx\dt\bigg]^p<\infty,
\end{align*}
for all $1\leq p<\infty$. Note that the last integrand is a lower bound for $H(\rho,\overline\rho)$ due to the elementary inequality
\begin{equation}
\begin{aligned}
\label{HrRho}
H(\rho,r) 
\geq c(r)\,
\left\{
\begin{array}{ll}
\vert \rho- r\vert^2 &:\text{ if }r/2 \leq \rho\leq 2r, 
\\
1+\rho^\gamma &:\text{ else },
\end{array}
\right.
\end{aligned}
\end{equation}
which holds for any $r>0$.
\end{Remark}
\begin{Remark}
Although Theorem \ref{thm:Romeo1} was shown in \cite{Romeo1} for $n=3$, it also applies in general dimensions replacing the bound $\gamma>\frac{3}{2}$ by $\gamma>\frac{n}{2}$. 
\end{Remark}

\subsection{Inviscid limit}
\label{sec:dataNStoEuler}
In this section we give the main result concerning the relation of the systems
\eqref{P1}--\eqref{P2} and \eqref{P1NS}--\eqref{P3NS}. From a formal point of view, the Navier--Stokes system \eqref{P1NS}--\eqref{P3NS} converges to the Euler system \eqref{P1}--\eqref{P2} if $\nu,\lambda\rightarrow0$. In order to make this idea rigorous, we have to analyse a singular limit.
Singular limit arguments for analysing the interactions between fluid dynamic models arises from suitable change of variables in time and space or by using dimensional analysis. Such transformations are now standard and interested readers can refer to \cite{Alazard2006low, lions1996mathematical} and the references within for further information.
\\
To study the inviscid limit result for \eqref{P1NS}--\eqref{P3NS}, we are interested in the transformation that leads to the following mappings:
\begin{align*}
\varrho \mapsto \rho_\varepsilon,\quad\mathbf{u} \mapsto\mathbf{v}_\varepsilon, \quad
 \nu \mapsto \varepsilon\nu,\quad \lambda\mapsto\varepsilon\lambda.
\end{align*}
This yields  the system
\begin{equation}
\begin{aligned}
\label{comprSPDE0}
\mathrm{d}\rho_\varepsilon + \mathrm{div}(\rho_\varepsilon \mathbf{v}_\varepsilon)\mathrm{d}t &= 0,
 \\
\mathrm{d}(\rho_\varepsilon \mathbf{v}_\varepsilon) + \Big[\mathrm{div}(\rho_\varepsilon  \mathbf{v}_\varepsilon\otimes \mathbf{v}_\varepsilon)     + a\nabla \rho^\gamma_\varepsilon \Big]\mathrm{d}t
& = 
\varepsilon\,\mathrm{div} \,\mathbb{S}(\nabla\mathbf{v}_\varepsilon) \,\mathrm{d}t 
+
 \mathbb{G}(\rho_\varepsilon,\rho_\varepsilon \mathbf{v}_\varepsilon)\mathrm{d}W,
\end{aligned}
\end{equation}
where the parameter $\varepsilon\in (0,1]$ corresponds to the inverse of the \textit{Reynolds number}. Our aim is to pass to the limit $\varepsilon\rightarrow0$.
We consider the following ill-prepared data that connects the inputs of Navier--Stokes and Euler system. We assume that the initial data
$(\rho_{0,\varepsilon},\mathbf{v}_{0,\varepsilon})$
 of the system \eqref{comprSPDE0} 
 satisfy the following conditions
\begin{equation}
\begin{aligned}
\label{intialCNSdata1}
\E\int_{\R^n}H(\rho_{0,\ep},\overline\varrho)\dx<\infty,\quad \rho_{0,\ep}|\bfv_{0,\ep}|^2\in L^{1}(\R^n),\quad 0< \varrho^-\leq \rho_{0,\ep}\leq\,\varrho^+\quad\mathbb P\text{-a.s.},
\end{aligned}
\end{equation}
where $\varrho^{-}$ and $\varrho^+$ are independent of $\ep$.
The initial data $(\varrho_0, \mathbf{u}_0)$ of the limit system \eqref{P1}--\eqref{P2} satisfy
\begin{equation}
\begin{aligned}
\label{intialCNSdata2}
(\varrho_0, \mathbf{u}_0)\in \overline\varrho+W^{s,2}(\mt)\times W^{s,2}(\mt),\quad \varrho_0\geq\varrho^->0\quad\mathbb P\text{-a.s.}
\end{aligned}
\end{equation}
Finally, we suppose that
\begin{equation}
\label{intialCNSdata3}
\begin{aligned}
\E\int_{\R^n}H(\rho_{0,\ep},\varrho_{0})\dx\xrightarrow{\varepsilon\searrow 0} 0,\quad \E\int_{\R^n}|\bfv_{0,\ep}-\bfu_0|^2\dx\xrightarrow{\varepsilon\searrow 0} 0.
\end{aligned}
\end{equation}
%
Our main result reads as follows.

\begin{Theorem}
\label{thm:invisicdLimit}
Let $\overline\varrho>0$ be given and suppose that \eqref{FG1}--\eqref{FG3} hold with $s>\frac{n}{2}+2$. Let $(\Omega,\mathfrak{F},\mathbb{P})$ be a complete probability space and $W$ a cylindrical Wiener process on $(\Omega,\mathfrak{F},\mathbb{P})$. Assume that
\begin{align}
\label{navierweaksol}
\big[(\Omega,\mathfrak{F},(\mathfrak{F}_t)_{t\geq0},\mathbb{P}),\rho_\varepsilon,\mathbf{v}_\varepsilon,W\big]_{\varepsilon>0}
\end{align}
is a family of finite energy weak martingale solution to the system \eqref{comprSPDE0} in the sense of Definition \ref{def:weakMartin} with $\overline\rho=\overline\varrho>0$.
On the same stochastic basis $(\Omega,\mathfrak{F},(\mathfrak{F}_t),\mathbb{P})$, 
consider the unique maximal strong pathwise solution to the Euler system \eqref{P1}--\eqref{P2} given by $(\varrho,\mathbf{u},(\mathfrak{t}_R)_{R\in\mathbb{N}},\mathfrak{t})$ driven by the same cylindrical Wiener process $W$. Assume that the initial data $(\rho_{0,\ep},\bfv_{0,\ep})$ and 
$(\varrho_0,\bfu_0)$ are $\mathfrak F_0$-measurable and satisfies  \eqref{intialCNSdata1}--\eqref{intialCNSdata3}.
Then we have
\begin{equation}
\label{functionalE}
\sup_{t\in(0,T)}\mathbb{E}\int_{\mathbb{R}^n}\bigg[ \frac{1}{2}\rho_\varepsilon \vert  \mathbf{v}_\varepsilon - \mathbf{u}\vert^2 + H(\rho_\varepsilon,\varrho)  \bigg](t\wedge \mathfrak{t}_R, \cdot)\,\mathrm{d}x\rightarrow 0
\end{equation}
as $\varepsilon\rightarrow0$ for all $R\in\mathbb N$.
\end{Theorem}
Before we proof Theorem \ref{thm:invisicdLimit}, we remark that \eqref{functionalE} implies that we have
\begin{align}
\rho_\varepsilon(\cdot\wedge\mathfrak{t}_R) \rightarrow \varrho(\cdot\wedge\mathfrak{t}_R) &\text{ in } L^{\overline\gamma}\big( \Omega\times(0,T);L^{\overline{\gamma}}_{loc}(\mathbb{R}^n) \big), \label{rhoGamma}\\
(\rho_\varepsilon\mathbf{v}_\varepsilon)(\cdot\wedge\mathfrak{t}_R) \rightarrow (\varrho \mathbf{u})(\cdot\wedge\mathfrak{t}_R) &\text{ in } L^{\frac{2\overline{\gamma}}{\overline{\gamma}+1}}\big( \Omega\times (0,T);L^{\frac{2\overline{\gamma}}{\overline{\gamma}+1}}_{loc}(\mathbb{R}^n) \big), \label{momtum}
\end{align}
where $\overline{\gamma}=\min\{2,\gamma\}$. The convergence \eqref{rhoGamma} follows from inequality \eqref{HrRho}
whilst \eqref{momtum} follows from the identity
\begin{align*}
\rho_\varepsilon\mathbf{v}_\varepsilon -\varrho\mathbf{u}= (\rho_\varepsilon-\varrho)\mathbf{u}+ \sqrt{\rho_\varepsilon}\sqrt{\rho_\varepsilon}(\mathrm{v}_\varepsilon -\mathbf{u})
\end{align*}
and H\"older's inequality, cf.  \cite{Sueur}.

\section{Proof of Theorem \ref{thm:main}}
\subsection{Approximation}

On a formal level, it can be seen that \eqref{P1}--\eqref{P2} is equivalent to

\begin{equation} \label{E30}
\D r + \vu \cdot \Grad r \ \dt + \frac{\gamma - 1}{2} r \Div \vu\dt = 0,
\end{equation}
\begin{equation} \label{E40}
\D \vu + \left[ \vu \cdot \Grad \vu + r \Grad r \right]  \dt  =  \tn{F} (r, \vu) \D W,
\end{equation}
where
$$ \varrho(r)=\Big(\frac{\gamma-1}{2a\gamma}r^2\Big)^{\frac{1}{\gamma-1}}, \quad \tn{F} (r, \vu) = \frac{1}{\vr(r)} \tn{G} (\vr(r), \vr(r) \vu),$$
cf. \cite[Sec. 2.4]{BrFeHo2016}. This can be made rigorous as long as $r$ (or equivalently, $\varrho$) is strictly positive. As we will solve \eqref{P1}--\eqref{P2} with respect to far field conditions, we seek a solution $r\in \overline r+W^{s,2}(\mt)$, where $\overline r=\sqrt{ \frac{2 a \gamma}{\gamma - 1} }\overline\vr^{\frac{\gamma-1}{2}}$.
In order to do so,
 we set $\hat r:=r-\overline r$ and aim to solve
\begin{equation} \label{E3}
\D \hat r + \vu \cdot \Grad \hat r \ \dt + \frac{\gamma - 1}{2} (\hat r+\overline r) \Div \vu\dt = 0,
\end{equation}
\begin{equation} \label{E4}
\D \vu + \left[ \vu \cdot \Grad \vu + (\hat r+\overline r) \Grad \hat r \right]  \dt  =  \hat{\tn{F}} (\hat r, \vu) \D W,
\end{equation}
where $\hat{\tn{F}} (\hat r, \vu)=\tn{F} (\hat r+\overline r, \vu)$.
We remark that the left-hand side of \eqref{E4} corresponds to a symmetric hyperbolic system, cf. Majda \cite{Majd}. In the stochastic case such system have been studied
in \cite{JUKim}. Unfortunately, the result from \cite{JUKim} does not apply to the general assumptions on $\tn{G}$ we have in mind.
In fact, the assumptions on the noise coefficient $\tn{F}$ are violated for small (close to zero) or large values of $r$. Due to this we replace $\tn{F}$ by
$$\tn{F}_R (r, \vu) = \frac{1}{\vr(r)} \varphi_R(\vr(r))\varphi_R(\vr(r)^{-1})\tn{G} (\vr(r), \vr(r) \vu),$$
where $\varphi_R:[0,\infty)\rightarrow[0,1]$ are smooth cut-off functions satisfying
\begin{equation}
\begin{aligned}
\label{cutOff}
\varphi_R(y)=\begin{cases}1,\quad &0\leq y\leq R,\\
0,\quad & R+1\leq y,
\end{cases}
\end{aligned}
\end{equation}
and similarly $\hat{\tn{F}}_R (\hat r, \vu) =\tn{F}_R (\hat r+\overline r, \vu) $.
We now study the system
\begin{equation} \label{E3R}
\D \hat r + \vu \cdot \Grad\hat r \ \dt + \frac{\gamma - 1}{2} (\hat r+\overline r) \Div \vu\dt = 0,
\end{equation}
\begin{equation} \label{E4R}
\D \vu + \left[ \vu \cdot \Grad \vu + (\hat r+\overline r) \Grad \hat r \right]  \dt  =  \hat{\tn{F}}_R (\hat r, \vu) \D W.
\end{equation}
We assume for the moment that
\begin{equation}\label{eq:vr0}
\|(\vr_0,\bfu_0) \|_{W^{s,2}(\mt)} < M, \ \vr_0 > \frac{1}{M} \ \prst\mbox{-a.s.}
\end{equation}
for some deterministic constant $M>0$. These assumptions will be relaxed later.
By definition of $\tn{F}_R$, the noise disappears if $\varrho(r)$ is larger than $R+1$ or smaller than $\frac{1}{R+1}$. Consequently, \eqref{FG1} and \eqref{FG2} imply that $\tn{F}_R$ (and $\hat{\tn{F}}_R$) is globally Lipschitz continuous on $W^{s,2}(\mt)\times W^{s,2}(\mt)$ for any fixed $R$. By \cite[Thm. 1.2]{JUKim} (where $u$ takes the role of $(\hat r,\bfu)$), there is a unique strong solution $(\hat r_R,\bfu_R,\mathfrak{t}_R)$ to \eqref{E3R}--\eqref{E4R} in the following sense:\footnote{The result from \cite{JUKim} requires the assumption $s>\frac{n}{2}+2$.}
\begin{itemize}
\item[(i)] $(\hat r_R,\bfu_R)$ is a $W^{s,2}(\mt)\times W^{s,2}(\mt)$-valued right-continuous progressively measurable process;
\item[(ii)] $\mathfrak{t}_R$ is a stopping time with respect to $(\mathfrak{F}_t)$ such that $\prst$-a.s.
\begin{align}\label{Kim1.8}
\mathfrak{t}_R=\lim_{N\rightarrow\infty} \mathfrak{t}_{R,N}
\end{align}
where 
\begin{align}\label{Kim1.9}
\mathfrak{t}_{R,N}=\inf\big\{0\leq t<\infty:\,\|(\hat r_R,\bfu_R)(t)\|_{s,2}\geq N\big\}
\end{align}
with the convention that $\mathfrak{t}_{R,N}=\infty$ if the set above is empty;
\item[(iii)] there holds $\prst$-a.s.
$$ (\hat r_R,\vu_R)(\cdot\wedge \mathfrak{t}_R) \in   C([0,T]; W^{s,2}(\tor)\times  W^{s,2}(\tor))$$
as well as
\[
\begin{split}
\hat r(t\wedge \mathfrak{t}_{R,N})  &=r_0-\int_0^{t \wedge \mathfrak{t}_{R,N}} \vu \cdot \Grad \hat r \ \dif s- \int_0^{t \wedge \mathfrak{t}_{R,N}}\frac{\gamma - 1}{2} (\hat r+\overline r) \Div \vu\ \dif s ,\\
\vu(t\wedge \mathfrak{t}_{R,N}) &= \vu_0- \int_0^{t \wedge \mathfrak{t}_{R,N}} \vu \cdot \Grad \vu   \ \dif s - \int_0^{t \wedge \mathfrak{t}_{R,N}} (\hat r+\overline r) \Grad \hat r   \ \dif s  +  \int_0^{t \wedge \mathfrak{t}_{R,N}}\hat{\tn{F}}_R (\hat r, \vu) \D W,
\end{split}
\]
for all $t\in[0,T]$ and all $N\geq1$, where $r_0=\sqrt{\tfrac{2a\gamma}{\gamma-1}}\varrho_0^{\frac{\gamma-1}{2}}-\overline r$.
\end{itemize}
Given $(\hat r_R,\bfu_R)$ we see that $ (r_R,\bfu_R)=(\hat r_R+\overline r,\bfu_R)$ solves  $\mathbb P$-a.s. 
\begin{equation} \label{E3b}
\D r_R + \vu_R \cdot \Grad r_R \ \dt + \frac{\gamma - 1}{2} r_R \Div \vu_R\dt = 0,
\end{equation}
\begin{equation} \label{E4b}
\D \vu_R + \left[ \vu_R \cdot \Grad \vu_R + r_R \Grad r_R \right]  \dt  =  \tn{F}_R (r_R, \vu_R) \D W.
\end{equation}
The aim in the following is to pass to the limit  $R\rightarrow\infty$ in \eqref{E3b} and \eqref{E4b}. This will be done by verifying the assumptions
of Lemma \ref{lem:GH}. We start by showing certain a priori estimates as a consequence of which we obtain \eqref{GH2}. Eventually we show uniqueness of
\eqref{E30}--\eqref{E40} which implies \eqref{GH1}.

\subsection{A priori estimates}
We immediately see that $(\vr_R,\bfu_R):=\Big( \Big(\frac{\gamma-1}{2a\gamma}r_R^2\Big)^{\frac{1}{\gamma-1}},\bfu_R\Big)$ solves  $\mathbb P$-a.s.
\begin{equation} \label{P1R}
\D \vr_R + \Div (\vr_R \vu_R) \ \dt = 0
\end{equation}
with $\vr_R(0)=\vr_0$ and far field $\overline\varrho$.
The standard maximum principle (see, e.g., \cite[eq. (7)]{DL} or \cite[eq. (2.8)]{FNP}) applied to \eqref{P1R} yields
\begin{align}\label{max}
\inf_{\mt}\varrho_0\exp\bigg({-\int_0^{\mathfrak{t}_{R,N}}\|\Div\bfu_R\|_\infty\dif s}\Big)\leq \varrho_R(t,x)\leq\,\sup_{\mt} \varrho_0\,\exp\Big({\int_0^{\mathfrak{t}_{R,N}}\|\Div\bfu_R\|_\infty\dif s}\bigg)
\end{align}
$\mathbb P$-a.s. for all $(t,x)\in(0,\mathfrak{t}_{R,N})\times\mt$.
Consequently, the definition of $\mathfrak{t}_{R,N}$, the embedding
$W^{s,2}(\mt)\hookrightarrow C^{1,\alpha}(\mt)$ (for $s>\frac{n}{2}+1$ and some $\alpha>0$), as well as \eqref{eq:vr0} implies that
\begin{align}\label{eq:maxr}
c_N^{-1}\inf_{\mt}\varrho_0\leq \varrho_R(t,x)\leq\,c_N\sup_{\mt} \varrho_0
\end{align}
$\mathbb P$-a.s. for all $(t,x)\in(0,\mathfrak{t}_{R,N})\times\mt$ and for a positive constant $c_N=c_N(\overline \varrho,M)$ which is independent of $R$. For a given $N\in\mathbb N$ there is $R_N$ such that \eqref{eq:maxr} implies
that 
\begin{align}\label{eq:maxrb}
R_N^{-1}\leq\varrho_{R_N}\leq R_N
\end{align} 
$\mathbb P$-a.s. in $(0,\mathfrak{t}_{R_N,N})\times\mt$.
\\Let $\alpha$ be a multiindex such that $|\alpha|\leq s$. Differentiating \eqref{E3b} in the $x$-variable, we obtain
\begin{align} \label{E8'}
\begin{aligned}
\D \partial^\alpha_x r_R &+ \left[\vu_R \cdot \Grad \partial^\alpha_x r_R \ + \tfrac{\gamma - 1}{2} \,r_R\, \Div \partial^\alpha_x  \vu_R \right] \dt  \\
&=\big[ \vu_R \cdot \partial^\alpha_x \Grad  r_R  - \partial^\alpha_x \left( \vu_R \cdot \Grad  r_R \right) \big] \dt\\
&+ \tfrac{\gamma - 1}{2}\left[r_R \partial^\alpha_x \Div  \vu_R - \partial^\alpha_x \left( r_R  \Div  \vu_R \right) \right] \dt\\
&=:T_1^R \dt+T^R_2 \dt.
\end{aligned}
\end{align}
Similarly, we differentiate \eqref{E4b} and deduce
that
\begin{align} \label{E9'}
\begin{aligned}
\D  \partial^\alpha_x \vu_R  &+\left[ \vu_R \cdot \Grad \partial^\alpha_x \vu_R + r_R \Grad \partial^\alpha_x r_R \right]  \dt\\
  &=  \left[ \vu_R \cdot \partial^\alpha_x \Grad  \vu_R - \partial^\alpha_x \left( \vu_R \cdot \Grad  \vu_R \right)  \right]  \dt  \\
&+  \left[ r_R \partial^\alpha_x \Grad  r_R -
\partial^\alpha_x \left( r_R  \Grad  r_R \right) \right]  \dt  \\
& +   \partial^\alpha_x \tn{F}_R (r_R, \vu_R) \D W\\
&=:T_3^R\dt+T^R_4 \dt + \partial^\alpha_x \tn{F}_R (r_R, \vu_R) \D W.
\end{aligned}
\end{align}
It follows from \eqref{E6} that the ``error'' terms may be handled as
\begin{equation} \label{E10'}
\begin{split}
\left\| T_1^R\right\|_{2} & \lesssim  \Big[
\| \Grad \vu_R \|_{\infty} \| \Grad^s r_R \|_{2} + \left\| \Grad r_R \right\|_{\infty} \| \Grad^s \vu_R \|_{2} \Big], \\
\left\| T_2^R \right\|_{2} & \lesssim  \Big[
\| \Grad r_R \|_{\infty} \| \Grad^s \vu_R  \|_{2} + \left\| \Div \vu_R  \right\|_{\infty} \| \Grad^s r_R \|_{2} \Big],
\\
\left\|  T_3^R \right\|_{2} & \lesssim 
\| \Grad \vu_R \|_{\infty} \| \Grad^s \vu_R  \|_{2} , \\
\left\|  T_4^R \right\|_{2} & \lesssim
\| \Grad r_R \|_{\infty} \| \Grad^s r  _R\|_{2}.
\end{split}
\end{equation}
Multiplying \eqref{E8'} by $\partial^\alpha_x (r_R-\overline r)$, we observe
\begin{align} \label{E12'}
&\left\| \partial^\alpha_x\big( r_R (t)-\overline r\big) \right\|_{2}^2 + (\gamma - 1)\int_0^t \int_{\mt} r_R \Div \partial^\alpha_x \vu_R \partial^\alpha_x (r_R-\overline r) \dx\ds \\
&\leq \left\| \partial^\alpha_x \big(r_0-\overline r) \right\|_{2}^2 +c\int_0^t
\left( \| \nabla_x\vu_R \|_{ \infty} \| r_R-\overline r \|_{s,2} + \| \nabla_x (r_R-\overline r) \|_{\infty} \| \vu _R\|_{s,2} \right)
\|\partial^\alpha_x (r_R-\overline r)\|_2\,\ds\nonumber
\end{align}
provided $|\alpha| \leq s$. Here, we took into account
$$
\int_{\mt}\vu_R \cdot \Grad \partial^\alpha_x r_R \partial^\alpha_x (r_R-\overline r)\,\dx=-\frac{1}{2}\int_{\mt}\Div\vu_R|\partial^\alpha_x (r_R-\overline r)|^2\,\dx
$$
as well as \eqref{E10'}.
To apply the same treatment to \eqref{E9'}, we apply It\^o's formula to the function $\int_{\mt}|\partial^\alpha_x\bfu_R|^2\dx$ and gain
\begin{equation} \label{E13'}
\begin{split}
 \left\|\partial^\alpha_x \vu_R(t) \right\|^2_2 \dx &+2\int_0^t\int_{\mt} \left[ \vu_R \cdot \Grad \partial^\alpha_x \vu_R + r_R \Grad \partial^\alpha_x (r_R -\overline r) \right] \cdot \partial^\alpha_x \vu_R \dx\ds \\
&=  \left\| \partial^\alpha_x \vu_0 \right\|^2  + 2\int_0^t\int_{\mt}\left[ T_3^R+T_4^R\right] \cdot \partial^\alpha_x \vu_R  \dx\ds  \\
& + 2\int_0^t \int_{\mt} \partial^\alpha_x \tn{F}_R (r_R, \vu_R) \cdot \partial^\alpha_x \vu_R \ \D W\\
&+\sum_{k\geq1} \int_0^t \int_{\mt}| \partial^\alpha_x \vc{F}_{R,k} (r_R, \vu_R)|^2 \dx\ds.
\end{split}
\end{equation}
Integrating by parts yields
\[
\begin{split}
\int_{\mt} &\big[ \vu_R \cdot \Grad \partial^\alpha_x \vu_R + r _R \Grad \partial^\alpha_x (r _R -\overline r) \big] \cdot \partial^\alpha_x \vu_R \dx\\=& - \frac{1}{2} \int_{\mt} |\partial^\alpha_x \vu_R |^2 \Div \vu_R \dx
 - \int_{\mt}r_R \Div \partial^\alpha_x \vu_R \partial^\alpha_x (r_R-\overline r) \dx  \\&- \int_{\mt} \Grad (r_R-\overline r) \cdot \partial^\alpha \vu_R \partial^\alpha_x (r_R- \overline r)\dx
\end{split}
\]
Now we combine \eqref{E10'}--\eqref{E13'} (and multiply \eqref{E13'} by $\gamma-1$) and observe that the term containing $r_R\partial^\alpha_x (r_R-\overline r)\Div \partial^\alpha_x \vu_R$ on the left hand side cancels out. We therefore conclude that
\begin{align}\label{strong:mainest}
\begin{aligned}
&\left\| (r_R(t)-\overline r,\bfu_R(t)) \right\|^2_{s,2}  \\
&\quad \leq  \left\| (r_0 -\overline r,\bfu_0) \right\|^2_{s,2}+c
\int_0^t 
\| \nabla\vu_R \|_{\infty} \Big( \| r_R-\overline r \|^2_{s,2} + \| \vu_R \|^2_{s,2} \Big) \ds
 \\
&\quad +c
\int_0^t   \| \nabla (r_R-\overline r) \|_{\infty}\Big( \| r_R -\overline r\|_{s,2}^2+ \| \vu _R\|^2_{s,2} \Big) \ds
 \\&\quad +c\int_0^t\int_{\mt}  \partial^\alpha_x \tn{F}_R (r_R, \vu_R) \cdot \partial^\alpha_x \vu_R \dx\ \D W \\
&\quad +c\sum_{k\geq1}\int_0^t \int_{\mt} | \partial^\alpha_x \vc{F}_{R,k} (r_R, \vu_R)|^2 \dx \ds\\
& \quad =(I)_0+ (I)_1+\dots+ (I)_4.
\end{aligned}
\end{align}
We choose $R\geq R_N$, where $R_N$ is chosen in a way that the cut-offs in the definition of $\tn{F}_R$ are not seen for  $t\leq\mathfrak{t}_{R,N}$, recall \eqref{eq:maxrb}.
Now we take the supremum over $0\leq t\leq\mathfrak{t}_{R,N}\wedge\delta$, where $\delta>0$, and apply expectations.
Using the definition of $\mathfrak t_{R,N}$ we easily obtain
\begin{align*}
\E\sup_{0\leq t\leq\mathfrak{t}_{R,N}\wedge\delta}\big[(I)_1+(I)_2\big]\leq c(N)\,
\E \int_0^{\mathfrak{t}_{R,N}\wedge\delta} 
 \| (r_R-\overline r, \vu_R) \|^2_{s,2}\dt\leq\,c(N)\delta.
\end{align*}
In accordance with \eqref{eq:maxrb}, the cut-offs in the definition of 
$\tn{F}_R $ are not seen for $t\leq \mathfrak{t}_{R,N}$ if $R\geq R_N$. Moreover, the norms of $\partial^\alpha \tn{F}_R$ can be controlled by $N$. 
First of all, we have by \eqref{FG1}--\eqref{FG3}
\begin{align*}
\E\sup_{0\leq t\leq\mathfrak{t}_{R,N}\wedge\delta}(I)_4&= \E
\sum_{k\geq1}\int_0^{\mathfrak{t}_{R,N}\wedge\delta} \int_{\mathcal K}|\partial^\alpha_x\vc{F}_{R,k} (r_R, \vu_R)|^2\dx\ds\\
&\leq\,c\,\E \int_0^{\mathfrak{t}_{R,N}\wedge\delta} \int_{\mathcal K} (|\partial^\alpha_x r_R|^2+ |\partial_x^\alpha\vu_R|^2) \dx \ds\\
&\leq\,c\, \E \int_0^{\mathfrak{t}_{R,N}\wedge\delta} \int_{\mathcal K} (1+|\partial^\alpha_x (r_R-\overline r)|^2+ |\partial_x^\alpha\vu_R|^2) \dx \ds\\
&\leq\,c\, \E
\int_0^{\mathfrak{t}_{R,N}\wedge\delta} \big(1+
 \| (r_R-\overline r, \vu_R) \|^2_{s,2}\big)\ds\leq\,c(N)\delta.
\end{align*}
Using the Burkholder-Davis-Gundy inequality, we can estimate the stochastic integral in the same fashion. After applying expectations, we gain using \eqref{FG1}--\eqref{FG3}
\begin{align*}
\E\bigg[\sup_{0\leq t\leq\mathfrak{t}_{R,N}\wedge\delta}&\bigg|\int_0^t \int_{\mt} \partial^\alpha_x \tn{F} (r_R, \vu_R) \cdot \partial^\alpha_x \vu_R\,\dx \ \D W\bigg|\bigg]\\
&\lesssim
\E\bigg[\sum_{k\geq1}\int_0^{\mathfrak{t}_{R,N}\wedge\delta} \bigg(\int_{\mt} \partial^\alpha_x \vc{F}_k (r_R, \vu_R) \cdot \partial^\alpha_x \vu_R \dx\bigg)^2\dt\bigg]^{\frac{1}{2}}\\
&\lesssim
\E\bigg[\int_0^{\mathfrak{t}_{R,N}\wedge\delta} \bigg(\sum_{k\geq1}\|\partial^\alpha_x\vc{F}_k (r_R, \vu_R)\|^2_{s,2}\bigg) \|\vu_R\|^2_{s,2}\dt\bigg]^{\frac{1}{2}}\\
&\lesssim
\E\bigg[\int_0^{\mathfrak{t}_{R,N}\wedge\delta} \int_{\mathcal K}\big(|\partial_x^\alpha r_R|^2+|\partial_x^\alpha\bfu_R|^2\big)\dx \|\vu_R\|^2_{s,2}\dt\bigg]^{\frac{1}{2}}\\
&\lesssim
\E\bigg[\int_0^{\mathfrak{t}_{R,N}\wedge\delta} \int_{\mathcal K}\big(1+|\partial_x^\alpha (r_R-\overline r)|^2+|\partial_x^\alpha\bfu_R|^2\big)\dx \|\vu_R\|^2_{s,2}\dt\bigg]^{\frac{1}{2}}\\
&\lesssim
\E\bigg[\int_0^{\mathfrak{t}_{R,N}\wedge\delta} \big(1+\|r_R-\overline r\|^2_{s,2}+\|\vu_R\|^2_{s,2}\big) \|\vu_R\|^2_{s,2}\dt\bigg]^{\frac{1}{2}}\\
&\leq\,c(N)\sqrt{\delta}.
\end{align*}
Plugging all together, we have shown
\begin{align*}
&\E\sup_{0\leq t\leq\mathfrak{t}_{R,N}\wedge\delta}\left\| (r_R(t)-\overline r,\bfu_R(t)) \right\|^2_{s,2}  \leq \E \left\| (r_0 -\overline r,\bfu_0) \right\|^2_{s,2}+c(N)(\sqrt{\delta}+\delta).
\end{align*}
We obtain
\begin{align*}
\lim_{\delta\rightarrow0}\sup_{R}\mathbb P\bigg(\sup_{0\leq t\leq\mathfrak{t}_{R,N}\wedge\delta}\left\| (r_R(t)-\overline r,\bfu_R(t)) \right\|^2_{s,2}  > \left\| (r_0-\overline r ,\bfu_0) \right\|^2_{s,2}+1\bigg)=0
\end{align*}
which is \eqref{GH2}.

\subsection{Pathwise uniqueness}
\label{subsec:uniq}

We mimick the approach of the last subsection and have  a look at the difference of two solutions
$(r_R, \vu_R)$ and $(r_L, \vu_L)$ which satisfies
\begin{equation}\label{UU1}
\begin{split}
\D \partial^\alpha_x (r_R - r_L)  &= -  \partial^\alpha_x \left( \vu_R \cdot \Grad r_R  + \frac{\gamma-1}{2} r_R \Div \vc{u}_R \right) \dt
\\ &+  \partial^\alpha_x \left( \vu_L \cdot \Grad r_L  + \frac{\gamma-1}{2} r_L \Div \vc{u}_L \right) \dt,
\end{split}
\end{equation}
and
\[
\begin{split}
\D  \partial^\alpha_x (\vc{u}_R - \vc{u}_L )   &= - \partial^\alpha_x \Big( \vu_R \cdot \Grad \vu_R   + r_R \Grad r_R  \Big) \dt + \partial^\alpha_x \Big( \vu_L \cdot \Grad \vu_L   + r_L \Grad r_L  \Big) \dt\\
&+  \Big[  \partial^\alpha_x \mathbb{F}(r_R, \vc{u}_R) -  \partial^\alpha_x \mathbb{F}(r_L, \vc{u}_L)
\Big]  {\rm d}W.
\end{split}
\]
Multiplying (\ref{UU1}) by $\partial^\alpha_x (r_R - r_L)$, we get
\begin{equation}\label{UU2}
\begin{split}
\frac{1}{2} \D \left| \partial^\alpha_x (r_R - r_L) \right|^2  &= -  \partial^\alpha_x \left( \vu_R \cdot \Grad r_R  + \frac{\gamma-1}{2} r_R \Div \vc{u}_R \right) \partial^\alpha_x (r_R - r_L)\dt
\\ &+  \partial^\alpha_x \left( \vu_L \cdot \Grad r_L  + \frac{\gamma-1}{2} r_L \Div \vc{u}_L \right) \partial^\alpha_x (r_R - r_L)\dt.
\end{split}
\end{equation}
Similarly, using It\^{o}'s product rule, we obtain
\begin{equation} \label{UU3}
\begin{split}
\frac{1}{2} & \D  \left| \partial^\alpha_x (\vc{u}_R - \vc{u}_L ) \right|^2 \\  &= - \partial^\alpha_x \Big( \vu_R \cdot \Grad \vu_R   + r_R \Grad r_R ) \Big) \cdot \partial^\alpha_x (\vc{u}_R - \vc{u}_L ) \dt \\
&+ \partial^\alpha_x \Big( \vu_L \cdot \Grad \vu_L   + r_L \Grad r_L \Big) \cdot \partial^\alpha_x (\vc{u}_R - \vc{u}_L )\dt\\
&+  \Big[ \partial^\alpha_x \mathbb{F}_R(r_R, \vc{u}_R) - \partial^\alpha_x \mathbb{F}_L(r_L, \vc{u}_L)
\Big] \cdot \partial^\alpha_x (\vc{u}_R - \vc{u}_L )  {\rm d}W\\
&+ \frac{1}{2}\sum_{k\geq1}  \Big( \partial^\alpha_x \vc{F}_{R,k}(r_R, \vc{u}_R) -  \partial^\alpha_x \vc{F}_{L,k}(r_L, \vc{u}_L)
\Big)^2 {\rm d}t.
\end{split}
\end{equation}
We sum \eqref{UU2} and \eqref{UU3}, integrate over the physical space, and perform the same estimates as in the previous section. Note that the highest order terms in \eqref{UU2} read
\[
\begin{split}
 &\int_{ \tor}\left( \vu_R \cdot \Grad \partial^\alpha_x r_R - \vu_L \cdot \Grad \partial^\alpha_x r_L \right)
\partial^\alpha_x \left( r_R - r_L \right) \ \dx\\
+ \frac{\gamma - 1}{2}  &\int_{ \tor} \left( r_R \Div \partial^\alpha_{x} \vu_R - r_L \Div \partial^\alpha_{x} \vu_L \right) \partial^\alpha_x \left( r_R - r_L \right) \dx\\
=  &\int_{ \tor}\Big( (\vu_R - \vu_L) \cdot \Grad \partial^\alpha_x r_R
\partial^\alpha_x \left( r_R - r_L \right)
 + \frac{1}{2} \Div \vu_L
\left| \partial^\alpha_x \left( r_R - r_L \right) \right|^2 \Big)\ \dx\\
+ \frac{\gamma - 1}{2} &\int_{ \tor}  (r_R  - r_L) \Div \partial^\alpha_{x} \vu_L  \partial^\alpha_x \left( r_R - r_L \right) \dx \\
+ \frac{\gamma - 1}{2} &\int_{ \tor} r_R \Div \partial^\alpha_{x} (\vu_R - \vu_L)
\partial^\alpha_x \left( r_R - r_L \right) \dx.
\end{split}
\]
Here, the last integral
cancels out after  integration by parts, with its counterpart in \eqref{UU3}.
Summing over all $\alpha$ with $|\alpha|\leq s-1$ we deduce
\begin{align} 
\D &\left( \left\| r_R - r_L \right\|^2_{W^{s-1,2}} + \left\| \vu_R - \vu_L \right\|^2_{W^{s-1,2}} \right)\nonumber
\\
& \leq c(R,L)  \left[ \left( 1 +\| (r_R,\vu_R) \|_{W^{s,2}}^2+\| (r_L,\vu_L) \|_{W^{s,2}}^2\right)
\left( \left\| r_R - r_L \right\|^2_{W^{s-1,2}} + \left\| \vu_R - \vu_L \right\|^2_{W^{s-1,2}} \right) \right] \dt\nonumber\\
& +   \Big[  \partial^\alpha_x \mathbb{F}_R(r_R, \vc{u}_R) -  \partial^\alpha_x \mathbb{F}_L(r_L, \vc{u}_L)
\Big] \cdot \partial^\alpha_x (\vc{u}_R - \vc{u}_L )  {\rm d}W\nonumber\\
&+ \frac{1}{2}\sum_{k\geq1}  \Big( \partial^\alpha_x \vc{F}_{R,k}(r_R, \vc{u}_R) -  \partial^\alpha_x \vc{F}_{L,k}(r_L, \vc{u}_L)
\Big)^2 {\rm d}t,\label{UU4}
\end{align}
where $s > \frac{n}{2} + 2$. As the initial data coincide, we obtain by Gronwall's lemma for $R,L$ large enough
\begin{align*}
\E\sup_{0\leq t\leq \mathfrak{t}_{R,N}\wedge \mathfrak{t}_{L,N}}\Big(\left\| r_R - r_L \right\|^2_{W^{s-1,2}} + \left\| \vu_R - \vu_L \right\|^2_{W^{s-1,2}}\Big)=0.
\end{align*}
This certainly yields
\begin{align*}\E\sup_{0\leq t\leq \mathfrak{t}_{R,N}\wedge \mathfrak{t}_{L,N}}\Big(\left\| r_R - r_L \right\|^2_{W^{s,2}} + \left\| \vu_R - \vu_L \right\|^2_{W^{s,2}}\Big)=0.
\end{align*}
which implies \eqref{GH1}. 

\subsection{Conclusion}
As shown in the last two subsections, we are in the position to apply
Lemma \ref{lem:GH}. We infer the existence  of a $\mathbb P$-a.s. positive stopping time $\mathfrak t$ and a predictable process $(\mathfrak t,r,\bfu)$ such that
\begin{align}\label{eq:3101}
\E\sup_{0\leq t\leq \mathfrak{t}}\Big(\left\| r_R - r \right\|^2_{W^{s,2}} + \left\| \vu_R - \vu\right\|^2_{W^{s,2}}\Big)=0
\end{align}
as $R\rightarrow\infty$. By \eqref{eq:3101}, it is easy to pass to the limit in \eqref{E3b} and \eqref{E4b}. By It\^{o}'s formula, we conclude that $(\mathfrak t,\vr,\bfu):=\Big( \mathfrak t,\Big(\frac{\gamma-1}{2a\gamma}r^2\Big)^{\frac{1}{\gamma-1}},\bfu\Big)$
is a solution to \eqref{P1}--\eqref{P2} in the sense of Definition \ref{def:strsol}.\\
Assume that $[\varrho^1,\vu^1,\mathfrak{t}^1]$ and $[\varrho^2,\vu^2,\mathfrak{t}^2]$ are two local strong solutions with the same initial datum. Then
we obtain that $[\varrho^1,\vu^1]$ and $[\varrho^2,\vu^2]$ coincide a.s. as a direct consequence of \eqref{UU4} (note that the systems \eqref{P1}--\eqref{P2} and \eqref{E30}--\eqref{E40} are equivalent up to the stopping time). This also implies
that the blow-up time for two maximal strong solutions (in the sense of Definition \ref{def:maxsol}) coincide. So, maximal strong solutions are unique.
So far, we have assumed \eqref{eq:vr0} which is quite restrictive. This assumption
can be removed as in \cite[Sec. 4.3]{BrFeHo2016}. Finally, it is standard to extend the local strong solution to a maximal strong solution, cf. \cite{BMS, Elw, jac,MikRoz}. For our purposes, the method from \cite[Sec. 4.4]{BrFeHo2016} can be used. The proof of Theorem \ref{thm:main} is complete.

\section{Proof of Theorem \ref{thm:invisicdLimit}}
Let 
\begin{align*}
\big[(\Omega,\mathfrak{F},(\mathfrak{F}_t)_{t\geq0},\mathbb{P}),\rho_\varepsilon,\mathbf{v}_\varepsilon,W\big]_{\varepsilon>0}
\end{align*}
be a sequence of finite energy weak martingale solutions to \eqref{comprSPDE0}, existence of which is guaranteed by Theorem \ref{thm:Romeo1}.
Our aim is to pass to the limit $\varepsilon\rightarrow0$.

\subsection{Relative energy inequality}
The \emph{relative energy inequality} is a tool which enables
us to compare $\rho_\varepsilon,\mathbf{v}_\varepsilon$ with some smooth comparison functions. 
Let $\big(f \,,\,\mathbf{U}  \big)$ be a pair of stochastic processes which are adapted to the filtration $(\mathfrak{F}_t)_{t\geq0}$ and which satisfies
\begin{equation}
\begin{aligned}
\label{operatorBB}
\mathrm{d}f  &= D^d_tf\,\mathrm{d}t  + \mathbb{D}^s_tf\,\mathrm{d}W,
\\
\mathrm{d}\mathbf{U}  &= D^d_t\mathbf{U}\,\mathrm{d}t  + \mathbb{D}^s_t\mathbf{U}\,\mathrm{d}W.
\end{aligned}
\end{equation}
In the above, $D^d_tf$, $D^d_t\mathbf{U}$ are functions of $(\omega,t,x)$ and $\mathbb{D}^s_tf$, $\mathbb{D}^s_t\mathbf{U}$ belong to $L_2\big(\mathfrak{U};L^2(\mathbb{R}^n)  \big)$ for a.e $(\omega,t)\in\Omega\times [0,T]$. 
For the \emph{relative energy functional}
\begin{align*}
\mathcal{E}\big(\rho_\varepsilon,\mathbf{v}_\varepsilon\left\vert \right. f, \mathbf{U}  \big)=\int_{\mathbb{R}^n}\bigg[ \frac{1}{2}\rho_\varepsilon\vert  \mathbf{v}_\varepsilon - \mathbf{U}\vert^2 + H(\rho_\varepsilon,f)  \bigg]\,\mathrm{d}x,  
\end{align*}we have that for any $t\in(0,T)$,

\begin{equation}
\begin{aligned}
\label{relativeEntropy}
&\mathcal{E} \big(\rho_\varepsilon,\mathbf{v}_\varepsilon\left\vert \right. f, \mathbf{U}  \big)  
(t)
+\int_0^{t} \int_{\mathbb{R}^n}\varepsilon\,\big[\mathbb{S}(\nabla_x  \mathbf{v}_\varepsilon)- \mathbb{S}(\nabla_x  \mathbf{U})\big]: \left( \nabla_x  \mathbf{v}_\varepsilon - \nabla_x  \mathbf{U} \right)\,\mathrm{d}x  \,\mathrm{d}s
\\&\leq
\mathcal{E} \big(\rho_\varepsilon,\mathbf{v}_\varepsilon\left\vert \right. f, \mathbf{U}  \big)(0) +M_{RE}(t)  + \int_0^t\mathcal{R} \big(\rho_\varepsilon,\mathbf{v}_\varepsilon\left\vert \right. f, \mathbf{U}  \big)(s)\,\mathrm{d}s
\end{aligned}
\end{equation}
$\mathbb P$-a.s., where
\begin{equation}
\begin{aligned}
\label{remainderRE}
\mathcal{R} \big(\rho_\varepsilon,\mathbf{v}_\varepsilon\left\vert \right. f, \mathbf{U}  \big) 
&=
\int_{\mathbb{R}^n}\varepsilon\,\mathbb{S}(\nabla_x  \mathbf{U}): \left( \nabla_x  \mathbf{U} - \nabla_x  \mathbf{v}_\varepsilon \right)\,\mathrm{d}x
\\
&+
\int_{\mathbb{R}^n}\rho_\varepsilon\big(D^d_t \mathbf{U}+ \mathbf{v}_\varepsilon\cdot  \nabla_x  \mathbf{U}\big) \cdot \big(\mathbf{U}-\mathbf{v}_\varepsilon\big) \,\mathrm{d}x
\\
&+
\int_{\mathbb{R}^n}\big[(f-\rho_\varepsilon)P''(f)D^d_tf + \nabla_x  P'(f) \cdot (f\mathbf{U}-\rho_\varepsilon\mathbf{v}_\varepsilon)  \big] \,\mathrm{d}x
\\
&+
\int_{\mathbb{R}^n}\big[p(f)   - p(\rho_\varepsilon) \big]\mathrm{div}_x (\mathbf{U}) \,\mathrm{d}x
\\
&+
\frac{1}{2}
\sum_{k\in\mathbb{N}}
\int_{\mathbb{R}^n}\rho_\varepsilon\bigg\vert \frac{\mathbf{G}_k(\rho_\varepsilon,\rho_\varepsilon\mathbf{v}_\varepsilon)}{\rho_\varepsilon}  -\mathbb{D}^s_t\mathbf{U}(e_k)  \bigg\vert^2 \,\mathrm{d}x.
\end{aligned}
\end{equation}
Here, $M_{RE}$ is a real valued square integrable martingale
and $P(\varrho)=\frac{a}{\gamma-1}\varrho^\gamma$ is the pressure potential.
Let us finally specify  the appropriate smoothness assumptions
 we require for $D^d_tf$, $D^d_t\mathbf{U}$, $\mathbb{D}^s_tf$ and $\mathbb{D}^s_t\mathbf{U}$. We suppose that
\begin{align}\label{eq:smooth}
\begin{aligned}
(f-\overline\rho) \in C^\infty_c([0,T]\times\R^n), \ \vc{U} \in C^\infty_c([0,T]\times\R^n), \ \quad\text{$\mathbb{P}$-a.s.},\\
\E\bigg[\sup_{t \in [0,T] } \| f-\overline\rho \|_{W^{1,q}(\tor)}^2\bigg]^q  + \E\bigg[ \sup_{t \in [0,T] } \| \vc{U} \|_{W^{1,q}(\tor)}^2\bigg]^q \leq c(q) \quad\forall\,\, 2 \leq q < \infty,
\end{aligned}
\end{align}
\begin{equation} \label{bound}
0 < \underline{f} \leq f(t,x) \leq \overline{f} \quad\text{$\mathbb{P}$-a.s.}
\end{equation}
Moreover, $f$, $\vc{U}$ satisfy
\begin{align}
&D^d f, D^d \vc{U}\in L^q(\Omega;L^q(0,T;W^{1,q}(\mt))),\quad
\mathbb D^s f,\mathbb D^s \vc{U}\in L^2(\Omega;L^2(0,T;L_2(\mathfrak U;L^2(\tor)))),\nonumber\\
&\bigg(\sum_{k\geq 1}|\mathbb D^s f(e_k)|^q\bigg)^\frac{1}{q},\bigg(\sum_{k\geq 1}|\mathbb D^s \vc{U}(e_k)|^q\bigg)^\frac{1}{q}\in L^q(\Omega;L^q(0,T;L^{q}(\mt))).\label{new}
\end{align}
The \emph{relative energy inequality} \eqref{relativeEntropy} is a consequence of the \emph{energy inequality} \eqref{eq:energy}.
The proof relies on a sophisticated application of It\^{o}'s formula in infinite dimensions. The latter one can be found in \cite[Lemma 3.1]{BrFeHo2015A}
in case of periodic boundary conditions. For a corresponding statement in the current setting where the underlying domain is $\R^n$ we refer to \cite[Sec. 3.6]{Romeo2}.
In order to prove Theorem \ref{thm:invisicdLimit},
we choose
$(f,\mathbf U)=(\varrho(\cdot\wedge \mathfrak t_R),\bfu(\cdot\wedge \mathfrak t_R))$ where $(\varrho,\mathbf{u},(\mathfrak{t}_R)_{R\in\mathbb{N}},\mathfrak{t})$ is the unique maximal strong pathwise solution to \eqref{P1}--\eqref{P2} which exists by Theorem \ref{thm:main} (note that the $C^\infty_c$-assumption in \eqref{eq:smooth} can be relaxed by a standard approximation argument). Recall that the stopping time $\mathfrak t_R$ announces the blow-up and satisfies
\begin{equation*}
\sup_{t\in[0,\mathfrak{t}_R]}\|\vu(t)\|_{1,\infty}\geq R\quad \text{on}\quad [\mathfrak{t}<T] ;
\end{equation*}
Moreover, $(f,\mathbf U)=(\varrho,\bfu)$ satisfies an equation of the form \eqref{operatorBB}, where
\begin{align*}
D^d_tf  = -\mathrm{div}_x (\varrho\mathbf{u}),\quad \mathbb{D}^s_tf=0,
\quad
D^d_t\mathbf{U}  = -\mathbf{u}\cdot\nabla_x \mathbf{u} -\frac{1}{\varrho}\nabla_x  p(\varrho),
\quad
\mathbb{D}^s_t\mathbf{U}  = \frac{1}{\varrho}\mathbb{G} (\varrho,\varrho\mathbf{u}).
\end{align*}
By Theorem \ref{thm:main} and \eqref{FG1}--\eqref{FG3}, it is easy to see that
\eqref{bound} and \eqref{new} are satisfied for $t\leq \mathfrak t_R$. Note in particular the lower bound for $\varrho$ which follows from the maximum principle \eqref{max} and \eqref{intialCNSdata1}.
So, \eqref{relativeEntropy} holds and the remainder takes the form
\begin{align}
\mathcal{R} \big(\rho_\varepsilon,\mathbf{v}_\varepsilon\left\vert \right. \varrho, \mathbf{u}  \big) 
&=
\int_{\mathbb{R}^n}\varepsilon\,\mathbb{S}(\nabla_x  \mathbf{u}): \left( \nabla_x  \mathbf{u} - \nabla_x  \mathbf{v}_\varepsilon \right)\,\mathrm{d}x\nonumber
\\
&+
\int_{\mathbb{R}^n}\rho_\varepsilon\Big(-\mathbf{u}\cdot\nabla_x \mathbf{u} -\frac{1}{\varrho}\nabla_x  p(\varrho)+ \mathbf{v}_\varepsilon\cdot  \nabla_x  \mathbf{u}\Big) \cdot \Big(\mathbf{u}-\mathbf{v}_\varepsilon\Big) \,\mathrm{d}x\nonumber
\\
&+
\int_{\mathbb{R}^n}\big[-(\varrho-\rho_\varepsilon)P''(\varrho)\mathrm{div}_x (\varrho\mathbf{u})+ \nabla_x  P'(\varrho) \cdot (\varrho\mathbf{u}-\rho_\varepsilon\mathbf{v}_\varepsilon)  \big] \,\mathrm{d}x\nonumber
\\
&+
\int_{\mathbb{R}^n}\big[p(\varrho)-p'(\varrho)(\varrho-\rho_\varepsilon)   - p(\rho_\varepsilon) \big]\mathrm{div}_x (\mathbf{u}) \,\mathrm{d}x+\int_{\mathbb{R}^n}p'(\varrho)(\varrho-\rho_\varepsilon)\mathrm{div}_x (\mathbf{u}) \,\mathrm{d}x\nonumber
\\
&+
\frac{1}{2}
\sum_{k\in\mathbb{N}}
\int_{\mathbb{R}^n}\rho_\varepsilon\bigg\vert \frac{\mathbf{G}_k(\rho_\varepsilon,\rho_\varepsilon\mathbf{v}_\varepsilon)}{\rho_\varepsilon}  -\frac{\mathbf{G}_k(\varrho,\varrho\mathbf{u})}{\varrho}  \bigg\vert^2 \,\mathrm{d}x.\label{remainderRE'}
\end{align}
Note that we can write
$$\varrho-\rho_\varepsilon=(\varrho-\overline\varrho)-\big(\rho_\varepsilon-\overline{\varrho}\big)\mathbb I_{|\rho_\varepsilon-\overline{\varrho}|\leq 1}-\big(\rho_\varepsilon-\overline{\varrho}\big)\mathbb I_{|\rho_\varepsilon-\overline{\varrho}|> 1},$$
where we have $\p$-a.s.
\begin{align*}
\big(\rho_\varepsilon-\overline{\varrho}\big)\mathbb I_{|\rho_\varepsilon-\overline{\varrho}|\leq 1}&\in L^\infty(0,T;L^2(\R^n)),\\
\big(\rho_\varepsilon-\overline{\varrho}\big)\mathbb I_{|\rho_\varepsilon-\overline{\varrho}|> 1}&\in L^\infty(0,T;L^\gamma(\R^n)),
\end{align*}
see Remark \eqref{rem:new}. Consequently, all terms in \eqref{remainderRE'}
are well-defined due to the regularity of $(\varrho,\bfu)$.

\subsection{Estimating the remainder}
In order to estimate the remainder in \eqref{remainderRE'} we follow ideas from \cite[Sec. 4]{BrFeHo2015A}. We tacitly assume that
$t\leq \mathfrak t_R$ such that $\bfu,\nabla\bfu,\varrho$ and $\varrho^{-1}$ can be bounded in terms of $R$.  
By using the identities
\begin{align*}
\varrho\nabla_x  P'(\varrho) = \nabla_x  p(\varrho) ,
\quad
\varrho\partial_t P'(\varrho) = \partial_t p(\varrho) ,
\quad
-\partial_t \varrho =\mathrm{div}_x (\varrho\mathbf{u}),
\end{align*}
it holds that
\begin{align*}
\int_{\mathbb{R}^n}&\bigg[\frac{\rho_\varepsilon}{\varrho}\nabla_x  p(\varrho) \cdot (\mathbf{v}_\varepsilon - \mathbf{u})  - (\varrho -\rho_\varepsilon)P''(\varrho)\mathrm{div}_x (\varrho\mathbf{u}) + \nabla_x  P'(\varrho)\cdot (\varrho\mathbf{u}-\rho_\varepsilon\mathbf{v}_\varepsilon)  \bigg] \,\mathrm{d}x
\\
&=\int_{\mathbb{R}^n}\bigg[ \frac{\rho_\varepsilon}{\varrho}\nabla_x  p(\varrho)\cdot (\mathbf{v}_\varepsilon - \mathbf{u})  + (\varrho -\rho_\varepsilon)\partial_t \big[P'(\varrho)\big] + \nabla_x  p(\varrho)\cdot\mathbf{u} -  \frac{\rho_\varepsilon}{\varrho}\nabla_x  p(\varrho)\cdot\mathbf{v}_\varepsilon  \bigg] \,\mathrm{d}x
\\
&=\int_{\mathbb{R}^n}\bigg[ \partial_t p(\varrho) -\frac{\rho_\varepsilon}{\varrho} \partial_t p(\varrho) + \nabla_x  p(\varrho)\cdot \mathbf{u} -  \frac{\rho_\varepsilon}{\varrho}\nabla_x  p(\varrho)\cdot \mathbf{u}  \bigg] \,\mathrm{d}x
\\
&=\int_{\mathbb{R}^n} \left(\frac{\varrho-\rho_\varepsilon}{\varrho} \right) \big(\partial_t p(\varrho)  + \nabla_x  p(\varrho)\cdot\mathbf{u}\big) \,\mathrm{d}x.
\end{align*}
However, since $(\varrho,\bfu)$ is a strong solution to the continuity equation, it satisfies the \textit{strong} renormalized continuity equation
\begin{align*}
\partial_tp(\varrho) + \nabla_x  p(\varrho)\cdot\mathbf{u}=-\gamma p(\varrho)\mathrm{div}_x (\mathbf{u}).
\end{align*}
By combining this with the identity $\varrho\, p'(\varrho)=\gamma\, p(\varrho)$ yields
\begin{align*}
\int_{\mathbb{R}^n} \left(\frac{\varrho-\rho_\varepsilon}{\varrho} \right) \big(\partial_t p(\varrho)  + \nabla_x  p(\varrho)&\cdot\mathbf{u}\big) \,\mathrm{d}x
\leq
\int_{\mathbb{R}^n} (\rho_\varepsilon -\varrho )p'(\varrho)\,\mathrm{div}_x (\mathbf{u})\,\mathrm{d}x
\end{align*}
since $1/\gamma<1$. 
So, by collecting the above estimates, we can now deduce from \eqref{remainderRE'} that for each $R\in\mathbb{N}$,
\begin{equation}
\begin{aligned}
\label{relativeEntropy1}
&\mathcal{E} \big(\rho_\varepsilon,\mathbf{v}_\varepsilon\left\vert \right. \varrho, \mathbf{u}  \big)  
(t \wedge \mathfrak{t}_R)
+\varepsilon\, \int_0^{t \wedge \mathfrak{t}_R} \int_{\mathbb{R}^n}\Big(\mathbb{S}(\nabla_x  \mathbf{v}_\varepsilon)-\mathbb{S}(\nabla_x  \mathbf{u})\Big) : (\nabla_x  \mathbf{v}_\varepsilon-\nabla_x  \mathbf{u})\,\mathrm{d}x  \,\mathrm{d}s
\\&\leq
\mathcal{E} \big(\rho_\varepsilon,\mathbf{v}_\varepsilon\left\vert \right. \varrho, \mathbf{u}  \big)(0) + M_{RE}(t \wedge \mathfrak{t}_R)  + \int_0^{t \wedge \mathfrak{t}_R} \tilde{\mathcal{R}} \big(\rho_\varepsilon,\mathbf{v}_\varepsilon\left\vert \right. \varrho, \mathbf{u}  \big)(s)\,\mathrm{d}s,
\end{aligned}
\end{equation}
where now
\begin{equation}
\begin{aligned}
\label{remainderRE1}
\tilde{\mathcal{R}} \big(\rho_\varepsilon,\mathbf{v}_\varepsilon\left\vert \right. \varrho, \mathbf{u}  \big) 
=\varepsilon&\,  \int_{\mathbb{R}^n}\mathbb{S}(\nabla_x  \mathbf{u}): (\nabla_x  \mathbf{u}-\nabla_x  \mathbf{v}_\varepsilon)\,\mathrm{d}x  \\
&+
\int_{\mathbb{R}^n}\rho_\varepsilon\big(\mathbf{v}_\varepsilon - \mathbf{u}\big)\cdot  \nabla_x  \mathbf{u}\cdot \big(\mathbf{u}-\mathbf{v}_\varepsilon\big) \,\mathrm{d}x
\\
&-
\int_{\mathbb{R}^n}\big[p(\rho_\varepsilon) -(\rho_\varepsilon - \varrho)p'(\varrho)   - p(\varrho) \big]\mathrm{div}_x (\mathbf{u}) \,\mathrm{d}x
\\
&+
\frac{1}{2}
\sum_{k\in\mathbb{N}}
\int_{\mathbb{R}^n}\rho_\varepsilon\bigg\vert \frac{\mathbf{G}_k(\rho_\varepsilon,\rho_\varepsilon\mathbf{v}_\varepsilon)}{\rho_\varepsilon}  -\frac{\mathbf{G}_k(\varrho,\varrho\mathbf{u})}{\varrho}  \bigg\vert^2 \,\mathrm{d}x.
\end{aligned}
\end{equation}
Now we observe that 
\begin{equation}
\begin{aligned}
\label{r1est}
\bigg\vert \int_{\mathbb{R}^n} \rho_\varepsilon(\mathbf{v}_\varepsilon-\mathbf{u})\cdot \nabla_x \mathbf{u}\cdot \big(  \mathbf{u} -  \mathbf{v}_\varepsilon \big)\,\mathrm{d}x \bigg\vert
&\leq
\Vert \nabla_x \mathbf{u}\Vert_{L^\infty_x}
\int_{\mathbb{R}^n} \rho_\varepsilon \vert \mathbf{v}_\varepsilon -  \mathbf{u} \vert^2\,\mathrm{d}x\\
&\leq\,c(R)\,\mathcal{E} \big(\rho_\varepsilon ,\mathbf{v}_\varepsilon\left\vert \right. \varrho, \mathbf{u}  \big) 
\end{aligned}
\end{equation}
$\mathbb P$-a.s.
Since the identity
\begin{align*}
&p(\rho_\varepsilon) -(\rho_\varepsilon - \varrho)p'(\varrho)   - p(\varrho) 
=
(\gamma-1)H(\rho_\varepsilon,\varrho)
\end{align*}
holds, it follows that
\begin{equation}
\begin{aligned}
\label{r2est}
\bigg\vert \int_{\mathbb{R}^n} &\big[p(\rho_\varepsilon) -(\rho_\varepsilon-\varrho)p'(\varrho)   - p(\varrho) \big]\mathrm{div}_x (\mathbf{u}) \,\mathrm{d}x  \bigg\vert
\\
&\leq \,c\,
\Vert \mathrm{div}_x \,\mathbf{u}\Vert_{L^\infty_x}
 \int_{\mathbb{R}^n}H(\rho_\varepsilon,\varrho)\,\mathrm{d}x\\
&\leq\,c(R)\,\mathcal{E} \big(\rho_\varepsilon ,\mathbf{v}_\varepsilon\left\vert \right. \varrho, \mathbf{u}  \big).
\end{aligned}
\end{equation}
$\mathbb P$-a.s.
Similarly, by Young's inequality for bilinear forms
\begin{align} 
&\left|  \int_{\R^n} \ep\tn{S}(\Grad \vc{u}):(\Grad \vc{v}_\varepsilon-\Grad \vu) \dx \right| \nonumber\\ &\leq
\frac{1}{2}\ep  \int_{\R^n} \Big( \tn{S} (\Grad \vc{v}_\varepsilon) - \tn{S} (\Grad \vu) \Big): \big( \Grad \vc{v}_\varepsilon - \Grad \vu \big) \dx   + c\, \varepsilon \int_{\R^n} \left| \mathbb{S}(\Grad \vc{u}) \right|^2  \dx \nonumber\\ &\leq
\frac{1}{2}\ep \ \int_{\R^n} \Big( \tn{S} (\Grad \vc{v}_\varepsilon) - \tn{S} (\Grad \vu) \Big): \big( \Grad \vc{v}_\varepsilon - \Grad \vu \big)  \dx \dt + c(R)\varepsilon.\label{Step2}
\end{align}
$\mathbb P$-a.s.
Lastly, we rewrite
\begin{equation}
\begin{aligned}
\label{nio0}
&
\frac{1}{2}
\sum_{k\in\mathbb{N}}
\int_{\mathbb{R}^n}\rho_\varepsilon\bigg\vert \frac{\mathbf{G}_k(\rho_\varepsilon,\rho_\varepsilon\mathbf{v}_\varepsilon)}{\rho_\varepsilon}  -\frac{\mathbf{G}_k(\varrho,\varrho\mathbf{u})}{\varrho}  \bigg\vert^2 \,\mathrm{d}x  
\\
&=
\frac{1}{2}
\sum_{k\in\mathbb{N}}
\int_{\mathcal{K}} \chi_{\{\rho_\varepsilon \leq \varrho/2\}}\, \rho_\varepsilon \bigg\vert \frac{\mathbf{G}_k(\rho_\varepsilon,\rho_\varepsilon\mathbf{v}_\varepsilon)}{\rho_\varepsilon}  -\frac{\mathbf{G}_k(\varrho,\varrho\mathbf{u})}{\varrho}  \bigg\vert^2 \,\mathrm{d}x
\\
&+
\frac{1}{2}
\sum_{k\in\mathbb{N}}
\int_{\mathcal{K}} \chi_{\{ \varrho/2 < \rho_\varepsilon < 2\varrho\}}\, \rho_\varepsilon \bigg\vert \frac{\mathbf{G}_k(\rho_\varepsilon ,\rho_\varepsilon\mathbf{v}_\varepsilon)}{\rho_\varepsilon}  -\frac{\mathbf{G}_k(\varrho,\varrho\mathbf{u})}{\varrho}  \bigg\vert^2 \,\mathrm{d}x
\\
&+
\frac{1}{2}
\sum_{k\in\mathbb{N}}
\int_{\mathcal{K}} \chi_{\{  \rho_\varepsilon \geq 2\varrho\}}\, \rho_\varepsilon \bigg\vert \frac{\mathbf{G}_k(\rho_\varepsilon,\rho_\varepsilon\mathbf{v}_\varepsilon)}{\rho_\varepsilon}  -\frac{\mathbf{G}_k(\varrho,\varrho\mathbf{u})}{\varrho}  \bigg\vert^2 \,\mathrm{d}x
\\
&=: I_1+I_2+I_3
\end{aligned}
\end{equation}
for $\mathcal{K}\Subset\mathbb{R}^3$ (recall \eqref{FG3}). We can now use the inequality $\rho\leq 1+ \rho^\gamma$ and \eqref{HrRho} to conclude that
\begin{equation}
\begin{aligned}
\label{nio1}
I_1 &\leq 
c\,
\sum_{k\in\mathbb{N}}
\int_{\mathcal{K}} \chi_{\{\rho_\varepsilon \leq \varrho/2\}}\,\bigg( \frac{1}{\rho_\varepsilon}\,\vert \mathbf{G}_k(\rho_\varepsilon ,\rho_\varepsilon \mathbf{v}_\varepsilon) \vert^2  +\frac{\rho_\varepsilon}{\varrho^2}\,\vert \mathbf{G}_k(\varrho,\varrho\mathbf{u}) \vert^2  \bigg) \,\mathrm{d}x
\\
&\leq c\, \int_{\mathcal{K}} \chi_{\{\rho_\varepsilon \leq  \varrho/2\}}\big( \rho_\varepsilon + \rho_\varepsilon\vert  \mathbf{v}_\varepsilon\vert^2  + \rho_\varepsilon\vert  \mathbf{u}\vert^2\big) \,\mathrm{d}x
\\
&\leq c(R)\, \int_{\mathcal{K}} \chi_{\{\rho_\varepsilon \leq  \varrho/2\}}\big(1+ \rho^\gamma_\varepsilon +\rho_\varepsilon \vert  \mathbf{v}_\varepsilon -  \mathbf{u}\vert^2\big) \,\mathrm{d}x
\\
&\leq c(R)\,
\mathcal{E} \big(\rho_\varepsilon ,\mathbf{v}_\varepsilon\left\vert \right. \varrho, \mathbf{u}  \big)  .
\end{aligned}
\end{equation}
Similarly, we obtain by  \eqref{FG1}, \eqref{FG2}, \eqref{HrRho}, the bounds on $\varrho$ and the mean-value theorem
\begin{align*}
I_2&\leq\,\frac{1}{2}\,\sum_{k\geq1}\int_{\mt}\chi_{\frac{\varrho}{2}\leq \rho_\ep\leq 2\varrho}\rho_\ep\Big(\frac{\vc{G}_k(\rho_\ep,\rho_\ep\bfv_\ep)}{\rho_\ep}-\frac {\vc{G}_k(\varrho,\rho_\ep\bfv_\ep)}{\varrho}\Big)^2\dx\\&+\frac{1}{2}\,\sum_{k\geq1}\int_{\mt}\chi_{\frac{\varrho}{2}\leq \rho_\ep\leq 2\varrho}\rho_\ep\Big(\frac{\vc{G}_k(\varrho,\rho_\ep\bfv_\ep)}{\varrho}-\frac {\vc{G}_k(\varrho,\varrho\mathbf{u})}{\varrho}\Big)^2\dx\\
&\leq \,c(R)\,\int_{\mt}\chi_{\frac{\varrho}{2}\leq \rho_\ep\leq 2\varrho}\Big(|\rho_\ep-\varrho|^2(1+|\rho_\ep\bfv_\ep-\varrho\mathbf{u}|^2\Big)\dx\\
&\leq \,c(R)\,\int_{\mt}\chi_{\frac{\varrho}{2}\leq \rho_\ep\leq 2\varrho}\Big(|\rho_\ep-\varrho|^2(1+|\mathbf{u}|^2)+|\rho_\ep(\bfv_\ep-\mathbf{u})|^2\Big)\dx\\
&\leq\,c(R)\,\int_{\mt}\chi_{\frac{\varrho}{2}\leq \rho_\ep\leq 2\varrho}|\rho_\ep-\varrho|^2\dx+\,\int_{\mt}\rho_\ep|\bfv_\ep-\mathbf{u}|^2\dx\\
&\leq\,c(R)\,\mathcal E\Big(\rho_\ep,\bfv_\ep\Big|\varrho,\mathbf{u}\Big).
\end{align*}
Estimating $I_3$ in \eqref{nio0} is similar to \eqref{nio1}.
So, we can conclude from \eqref{nio0} that
\begin{equation}
\begin{aligned}
\label{r3est}
\frac{1}{2}
\sum_{k\in\mathbb{N}}
\int_{\mathbb{R}^n}\rho_\varepsilon \bigg\vert \frac{\mathbf{G}_k(\rho_\varepsilon ,\rho_\varepsilon \mathbf{v}_\varepsilon)}{\rho_\varepsilon}  -\frac{\mathbf{G}_k(\varrho,\varrho\mathbf{u})}{\varrho}  \bigg\vert^2 \,\mathrm{d}x  
\leq c(R)\,
\mathcal{E} \big(\rho_\varepsilon ,\mathbf{v}_\varepsilon\left\vert \right. \varrho, \mathbf{u}  \big).
\end{aligned}
\end{equation}

\subsection{Conclusion}
Collecting the estimates \eqref{r1est}--\eqref{r3est}, we have shown that
\begin{equation}
\begin{aligned}
\label{r4est}
\int_0^{t \wedge \mathfrak{t}_R}
 \mathcal{R} \big(\rho_\varepsilon ,\mathbf{v}_\varepsilon \left\vert \right. \varrho, \mathbf{u}  \big) \,\mathrm{d}s
\leq c(R)\bigg(\,\int_0^{t \wedge \mathfrak{t}_R}
\mathcal{E} \big(\rho_\varepsilon ,\mathbf{v}_\varepsilon \left\vert \right. \varrho, \mathbf{u}  \big)  
(s)\,\mathrm{d}s+\varepsilon\bigg).
\end{aligned}
\end{equation}
Combining \eqref{r4est} and \eqref{relativeEntropy1} and applying Gronwall's lemma yields
\begin{align}
\nonumber
\sup_{t\in(0,T)}\mathbb{E}\,  \Big[\mathcal{E} \big(\rho_\varepsilon ,\mathbf{v}_\varepsilon\left\vert \right. \varrho, \mathbf{u}  \big)  
(t \wedge \mathfrak{t}_R)\Big]
&+\varepsilon\, \mathbb{E}\,\Big[\int_0^{T \wedge \mathfrak{t}_R} \int_{\mathbb{R}^n}\big[\mathbb{S}(\nabla_x \mathbf{v}_\varepsilon)-\mathbb{S}(\nabla_x \mathbf{u})\big] : (\nabla_x \mathbf{v}_\varepsilon  - \nabla_x \mathbf{u})\,\mathrm{d}x  \,\mathrm{d}s \Big]
\\&\leq c(R)\,
\mathbb{E}\,\Big[\mathcal{E} \big(\rho_\varepsilon ,\mathbf{v}_\varepsilon\left\vert \right. \varrho, \mathbf{u}  \big)(0)+\varepsilon\Big].
\label{relativeEntropy2}
\end{align}
Note that we have
\begin{equation*}
\begin{aligned}
\mathcal{E} \big(\rho_\varepsilon ,\mathbf{v}_\varepsilon\left\vert \right. \varrho, \mathbf{u}  \big)  (0) 
&=
\int_{\mathbb{R}^n} \frac{1}{2}\rho_{0,\varepsilon}\big\vert \mathbf{v}_{0,\varepsilon} - \mathbf{u}_{0} \big\vert^2  \,\mathrm{d}x
+
\int_{\mathbb{R}^n}H\Big(\varepsilon\rho_{0,\varepsilon} ,\varrho_0 \Big) \,\mathrm{d}x
\end{aligned}
\end{equation*}
which converges to zero in expectation by \eqref{intialCNSdata1}--\eqref{intialCNSdata3}.
Consequently, we obtain
\begin{equation}
\begin{aligned}
\mathcal{E} \big(\rho_\varepsilon ,\mathbf{v}_\varepsilon\left\vert \right. \varrho, \mathbf{u}  \big)  (0) 
\rightarrow 0
\end{aligned}
\end{equation}
as $\varepsilon\rightarrow0$.
The convergence \eqref{functionalE} then follows from passing to the limit $\varepsilon\rightarrow0$ in \eqref{relativeEntropy2}.\hfill$\Box$

\end{document}